\documentclass[pdflatex,sn-mathphys-num]{sn-jnl}
\hypersetup{hidelinks,hypertexnames=false}
\makeatletter
\def\email#1{\global\advance\emailcnt by 1\relax%
\if@corauemail%
  \g@addto@macro\corrauthemail{\setcounter{footnote}{0}#1 }%
\else%
  \g@addto@macro\authemail{\setcounter{footnote}{0}#1;\ }%
\fi}
\makeatother

\usepackage{graphicx}%
\usepackage{multirow}%
\usepackage{amsmath,amssymb,amsfonts}%
\usepackage{amsthm}%
\usepackage{mathrsfs}%
\usepackage[title]{appendix}%
\usepackage{textcomp}%
\usepackage{manyfoot}%
\usepackage{booktabs}%
\usepackage{array}%
\usepackage{algorithm}%
\usepackage{algorithmicx}%
\usepackage{algpseudocode}%
\usepackage{listings}%

\usepackage{float}
\usepackage{url}
\usepackage{bm}

\newcommand{\eps}{\varepsilon}
\newcommand{\Om}{\Omega}

\newcommand{\R}{\mathbb{R}}
\newcommand{\dd}{\,\mathrm{d}}
\newcommand{\avg}[1]{\left\langle #1\right\rangle_Y}
\newcommand{\norm}[2][]{\left\lVert #2\right\rVert_{#1}}
\newcommand{\CD}[1]{{}^{C}\!D_t^{#1}}

\newcommand{\Ah}{\widehat A}
\newcommand{\rhoh}{\widehat\rho}
\newcommand{\muh}{\widehat\mu}
\newcommand{\Ly}{\mathcal L_{\by}}

\newcommand{\Ddisc}{\mathcal D}
\newcommand{\bx}{\bm{x}}
\newcommand{\by}{\bm{y}}

\newcommand{\bxi}{\bm{\xi}}

\theoremstyle{thmstyleone}%
\newtheorem{theorem}{Theorem}
\newtheorem{lemma}{Lemma}
\newtheorem{assumption}{Assumption}
\theoremstyle{thmstyletwo}%
\newtheorem{remark}{Remark}%

\theoremstyle{thmstylethree}%

\begin{document}

\title{A high-order multi-scale method and its error analysis for nonlocal-in-time fractionally damped wave problems in composite materials}


\author[1]{\fnm{Gehao} \sur{Wang}}

\author*[1]{\fnm{Hao} \sur{Dong}}\email{donghaoxd@xidian.edu.cn}

\author[2]{\fnm{Yufeng} \sur{Nie}}

\author[3]{\fnm{Junzhi} \sur{Cui}}

\affil[1]{\orgdiv{School of Mathematics and Statistics}, \orgname{Xidian University}, \orgaddress{\city{Xi'an}, \postcode{710071}, \country{PR China}}}

\affil[2]{\orgdiv{School of Mathematics and Statistics}, \orgname{Northwestern Polytechnical University}, \orgaddress{\city{Xi'an}, \postcode{710129}, \country{PR China}}}

\affil[3]{\orgdiv{LSEC, ICMSEC, Academy of Mathematics and Systems Science}, \orgname{Chinese Academy of Sciences}, \orgaddress{\city{Beijing}, \postcode{100190}, \country{PR China}}}


\abstract{Significant spatial heterogeneity and nonlocal history-dependent attenuation jointly govern wave propagation in composite materials, making full-scale direct numerical simulation (DNS) computationally prohibitive. This study proposed a high-order multi-scale (HOMS) method for nonlocal-in-time fractionally damped wave problems with highly oscillatory and discontinuous material coefficients. Using a multi-scale asymptotic expansion approach, the macroscopic homogenized equation, and first- and second-order correction terms are derived to construct the high-accuracy HOMS solution.
	Moreover, a local pointwise error analysis identifies the leading residual eliminated by this HOMS solution and validates its local physical balance preservation. Then, a global energy-norm error estimate establishes the \(O(\varepsilon)\)-order convergence rate for the HOMS solution. Subsequently, an efficient spatiotemporal multi-scale algorithm integrating the finite element and finite difference methods is designed, which employs a novel fully discrete \(L2\text{-}1_{\sigma}\)–Newmark scheme, and corresponding unconditional shifted-energy stability and explicit convergence rates are also obtained for the macroscopic homogenized displacement and velocity, respectively.
	Both two- and three-dimensional examples are presented to demonstrate the accuracy, efficiency, reduced cost, and robustness of the proposed HOMS method. This work provides an effective multi-scale tool for nonlocal-in-time wave simulation in engineering applications.}




\maketitle
\section{Introduction}
\label{sec:1}
As engineering applications demand effective vibration isolation and noise reduction, accurate prediction of the dynamic response of dissipative composite materials becomes important \cite{LjungMalqvistPersson2021}. These structures combine constituent materials with distinct properties, while their microscopic features are much smaller than the overall structural dimensions. Hence, periodic composite materials exhibit multi-scale spatial features, with rapidly varying material properties across the repeating microstructure \cite{CherednichenkoWaurick2018}. The present model also allows discontinuities across material interfaces. During wave propagation, dissipative constituents with temporal memory produce a response that depends on preceding motion. In some media, attenuation exhibits a non-integer power-law dependence on frequency, motivating a history-dependent description through fractional damping \cite{KaltenbacherRundell2022}. The interaction of spatial heterogeneity and temporal memory introduces additional difficulties in predicting local wave fields and their attenuation. DNS of the resulting nonlocal-in-time wave problem requires fine spatial meshes to resolve repeated material interfaces while accounting for the temporal history, leading to substantial computational demands. Therefore, efficient and accurate HOMS computational methods are important for simulating nonlocal-in-time wave propagation in dissipative composite materials.

For homogeneous media, established wave solvers include Newmark time integration for structural dynamics \cite{Newmark1959} and high-order implicit and implicit--explicit schemes on overset grids \cite{CarsonEtAl2025}. When damping depends on preceding motion, an instantaneous damping law must be supplemented by a history-dependent description. Fractional derivatives incorporate material memory into wave models through temporal convolution kernels. Mielke et al. derived a fractionally damped wave equation from a coupled wave--diffusion system and characterized its energy–dissipation structure \cite{MielkeNetzZendehroud2021}. Baker~\cite{BakerBanjai2022} and Cao~\cite{CaoAlMaskari2026} both established numerical error estimates for wave equations with fractional damping. Baker et al. also analyzed time stepping for the nonlinear Westervelt equation with time-fractional damping and corrections for initial singularities \cite{BakerBanjaiPtashnyk2024}. These studies provide an important analytical and numerical foundation for wave equations with fractional damping. However, applying these methods to composite materials still requires accounting for rapid spatial variations in material parameters and their effects on local wave fields and attenuation.

Let $\eps$ denote the microscopic periodicity scale. As $\eps\to0$, resolving material interfaces increases the spatial unknowns and the cost of storing their temporal histories. Fast convolution algorithms can reduce the temporal storage and computational burden \cite{ZengTurnerBurrage2018,BanjaiLopezFernandez2019}. These computational demands also motivate multi-scale methods that incorporate microstructural information into coarse-scale approximations. Representative approaches include the multiscale finite element method (MsFEM) \cite{HouWu1997}, the heterogeneous multiscale method (HMM) \cite{AbdulleEtAl2012}, the localized orthogonal decomposition method (LOD) \cite{MalqvistPeterseim2014}, and the FE$^2$ computational homogenization method \cite{FeyelChaboche2000}. These methods are generally somewhat less computationally efficient than homogeneous macroscopic models and less accurate than fully resolved microscopic models. When scale separation holds, the asymptotic homogenization method (AHM) determines effective homogenized material parameters through cell problems and scale averaging. The resulting macroscopic homogenized equation has spatially constant coefficients \cite{Allaire1992}, reducing the spatial unknowns involved in the history-dependent computation. However, the macroscopic homogenized solution alone generally cannot recover microscopic oscillations of local physical fields. Cui and co-workers developed HOMS methods with high-order corrections combining local auxiliary cell functions, the macroscopic homogenized solution, and its derivatives. These methods have been applied to heat conduction and mechanical problems in composite materials \cite{MaCuiLiWang2016,YangEtAl2017HOTS,DongEtAl2019Shell,DongEtAl2023NonlinearThermomechanical}.

In heterogeneous media, Dong and Cao developed high-order correctors and algorithms for second-order wave equations with rapidly oscillating coefficients \cite{DongCao2009}. Dong et al. constructed second-order two-scale expansions and a Newmark-based algorithm for damped waves in quasi-periodic composites \cite{DongEtAl2017DampedWave}. For strongly damped waves, Ljung et al. developed a generalized finite element method allowing independent spatial variations in damping and propagation coefficients \cite{LjungMalqvistPersson2021}. Subsequent studies address time-dependent multi-scale coefficients \cite{MaierVerfuerth2022}, fully explicit mass-lumped multi-scale discretization \cite{GeeversMaier2023}, and high-contrast media \cite{FressartVerfuerth2024}. Krumbiegel and Maier established energy-norm error estimates for high-order localized orthogonal decomposition combined with time integration \cite{KrumbiegelMaier2025}. Temporal memory must also be represented in the heterogeneous governing equation. Related analyses include operator-norm homogenization estimates for one-dimensional fractional elasticity \cite{CherednichenkoWaurick2018} and first-order approximation estimates for time-fractional diffusion with smooth periodic coefficients \cite{HuLi2020}. Bai et al. developed multicontinuum homogenization for a time-fractional diffusion--wave equation with high-contrast heterogeneous coefficients \cite{BaiEtAl2026}. For time discretization, the $L2$-$1_\sigma$ formula is second-order accurate for fractional diffusion under suitable regularity assumptions \cite{Alikhanov2015}. Second-order schemes also cover multi-term fractional diffusion--wave equations \cite{SunZhaoSun2019} and variable-coefficient fractional wave equations \cite{WangZheng2022}. To sum up, existing research on the HOMS method for modeling and computation of composite wave problems with classical inertia and Caputo fractional damping remains limited. Meanwhile, growing engineering demands for vibration isolation and noise reduction strongly motivate continued research on this challenging problem.

The central challenge is to reconstruct microscopic wave fields and their gradients from the macroscopic homogenized solution while accounting for spatial heterogeneity and temporal memory. Furthermore, in bounded domains, periodic cell correctors generally do not satisfy external boundary conditions, introducing a boundary-layer contribution to the global error. In addition, cell discretization errors must also be tracked through the homogenized material parameters to the macroscopic homogenized solution. The fully discrete scheme additionally requires stability and convergence analysis.

The principal contributions are threefold. First, we derive the macroscopic homogenized equation and first- and second-order auxiliary cell problems for the nonlocal-in-time scalar wave model. Second, the pointwise error analysis shows that the HOMS method cancels the leading residual of the LOMS method and preserves local physical balance. We also establish global energy-error estimates using homogeneous Dirichlet cell conditions. Third, the fully discrete $L2$-$1_\sigma$--Newmark scheme for the homogenized equation admits unconditional shifted-energy stability. Under the stated assumptions, displacement and velocity errors are $O(h_Y^{2p_Y}+h^p+\tau^2)$ in the $H^1$ and $L^2$ norms, respectively. Two- and three-dimensional examples demonstrate accurate reconstruction at reduced computational cost. A large-scale composite plate example further demonstrates feasibility when DNS exceeds available resources.

The remainder of this paper is organized as follows. Section~\ref{sec:2} introduces the governing problem and derives the macroscopic homogenized equation, auxiliary cell problems, and multi-scale asymptotic solutions. Section~\ref{sec:3} presents the local residual analysis and global convergence analysis. Section~\ref{sec:4} develops a high-precision finite element method coupled with an $L2$-$1_\sigma$--Newmark time-stepping scheme. Section~\ref{sec:5} establishes stability and convergence results. Section~\ref{sec:6} presents the numerical examples, and Section~\ref{sec:7} summarizes the conclusions.

\section{High-order multi-scale computational model for nonlocal-in-time wave problems}
\label{sec:2}
\subsection{Statement of the multi-scale problem}
\label{subsec:2.1}

Let $d\in\{2,3\}$, and let $\Om\subset\R^d$ be a bounded, connected
convex polygonal or polyhedral domain.  Within the framework of AHM, we introduce the microscopic coordinates
\(\boldsymbol y=\boldsymbol x/\varepsilon
=(x_1/\varepsilon,\ldots,x_d/\varepsilon)
=(y_1,\ldots,y_d),\)
and take $Y=[0,1]^d$ as the periodic unit cell (PUC). For a $Y$-periodic function $q$, write
$q^\eps(\bx)=q(\bx/\eps)$ and
\[
\avg{q}=|Y|^{-1}\int_Y q(\by)\dd \by.
\]
Here $|Y|:=\int_Y 1\,\dd\by$ denotes the $d$-dimensional measure of the PUC $Y$, and $|Y|^{-1}=1/|Y|$ is the normalization factor. For $Y=[0,1]^d$, we have $|Y|=1$. Throughout, we employ the Einstein summation convention.
We use $\partial_{x_k}:=\partial/\partial x_k$ and $\partial_{x_kx_\ell}:=\partial_{x_k}\partial_{x_\ell}$ for spatial derivatives, with analogous notation for higher orders and derivatives with respect to $\by$.
Subscripts $t$, $tt$, and their higher-order analogues denote time derivatives.

Let $X$ be a Banach space.  For a Bochner-integrable $z\in L^1(0,T;X)$ and every $\eta>0$, define
\begin{equation*}
	I^\eta z(t)=\frac{1}{\Gamma(\eta)}
	\int_0^t(t-s)^{\eta-1}z(s)\dd s.
\end{equation*}
Throughout the paper, $I^\eta$ denotes this integral operator, whereas $\Gamma(\eta)$ denotes Euler's gamma function.

For $1<\alpha<2$ and $v\in W^{2,1}(0,T;X)$, the Caputo derivative is
\begin{equation*}
	\CD{\alpha}v(t)
	=\frac{1}{\Gamma(2-\alpha)}\int_0^t
	(t-s)^{1-\alpha}v_{ss}(s)\dd s.
\end{equation*}
Set
\begin{equation*}
	\beta=\alpha-1\in(0,1),\qquad
	\gamma=2-\alpha=1-\beta\in(0,1).
\end{equation*}
For $0<\beta<1$ and $z\in W^{1,1}(0,T;X)$, define
\begin{equation*}
	\CD{\beta}z
	=\frac{\dd}{\dd t}I^{1-\beta}(z-z(0))
	=I^{1-\beta}z_t.
\end{equation*}
Consequently, if $v=u_t$, then
\begin{equation*}
	\CD{\alpha}u=I^{\gamma}u_{tt}=\CD{\beta}v,
	\qquad v\in W^{1,1}(0,T;X).
\end{equation*}

We consider the nonlocal‑in‑time fractionally damped wave problem:

\begin{equation}
	\label{eq:1}
	\left\{\begin{aligned}
		&\rho^\eps(\bx) u_{tt}^\eps(\bx,t)+\mu^\eps(\bx)\CD{\alpha}u^\eps(\bx,t) \\
		&\quad -\partial_{x_i}\!\left(a_{ij}^\eps(\bx)\partial_{x_j}u^\eps(\bx,t)\right)
		=f(\bx,t) \quad&&\text{in }\Om\times(0,T], \\
		&u^\eps(\bx,t)=\widehat u(\bx) &&\text{on }\partial\Om\times[0,T],\\
		&u^\eps(\bx,0)=g(\bx),\qquad u_t^\eps(\bx,0)=\psi(\bx) &&\text{in }\Om.
	\end{aligned}\right.
\end{equation}

 In model \eqref{eq:1}, $u^\varepsilon(\boldsymbol{x},t)$ denotes the unknown scalar displacement field.  The coefficients $\rho^\varepsilon(\boldsymbol{x})$, $\mu^\varepsilon(\boldsymbol{x})$, and $A^\varepsilon(\boldsymbol{x})=\{a_{ij}^\varepsilon(\boldsymbol{x})\}$ represent the mass density, fractional damping coefficient, and stiffness coefficient matrix, respectively. In terms of the microscopic coordinates $\boldsymbol{y}=\boldsymbol{x}/\varepsilon$, these coefficients can be expressed as
 \( \rho^\varepsilon(\boldsymbol{x})=\rho(\boldsymbol{y}),\ 
 \mu^\varepsilon(\boldsymbol{x})=\mu(\boldsymbol{y}),\ 
 A^\varepsilon(\boldsymbol{x})=A(\boldsymbol{y}),\)
 where $\rho$, $\mu$, and $A=\{a_{ij}\}$ are $1$-periodic functions of $\boldsymbol{y}$ and may be discontinuous across material interfaces. The term $\rho^\varepsilon u_{tt}^\varepsilon$ describes classical inertia, whereas $\mu^\varepsilon\,{}^C D_t^\alpha u^\varepsilon$ represents history-dependent damping, with $\alpha\in(1,2)$ denoting the fractional order. The prescribed source $f(\boldsymbol{x},t)$ represents the distributed force density. The functions $g(\boldsymbol{x})$ and $\psi(\boldsymbol{x})$ specify the initial displacement and velocity, respectively. The prescribed displacement $\widehat u(\bx)$  is imposed on the entire external boundary $\partial\Omega$. The parameter $T$ denotes the final time.

\begin{assumption}
	\label{ass:1}
	\par\smallskip\noindent\textup{(A1)}\quad
	The matrix $A=(a_{ij})$ is real and symmetric. There are positive constants $a_*,a^*$, $\rho_*,\rho^*$, and $\mu_*,\mu^*$. For almost every $\by\in Y$ and every $\bxi\in\R^d$, they satisfy
	\[
	\begin{aligned}
		a_*|\bxi|^2\le a_{ij}(\by)\xi_i\xi_j\le a^*|\bxi|^2,\\
		\rho_*\le\rho(\by)\le\rho^*,\qquad
		\mu_*\le\mu(\by)\le\mu^*.
	\end{aligned}
	\]
	\par\smallskip\noindent\textup{(A2)}\quad
	The coefficients $A$, $\rho$, and $\mu$ are $Y$-periodic and belong to	$L^\infty(Y)^{d\times d}$, $L^\infty(Y)$, and $L^\infty(Y)$, respectively.  The material partition is allowed to be piecewise smooth, and the coefficients may jump across its interfaces.	\par\smallskip\noindent\textup{(A3)}\quad The  boundary data satisfy $\widehat u\in H^{1/2}(\partial\Om)$. The initial displacement satisfies $g\in H^1(\Om)$ and $\operatorname{Tr}g=\widehat u$, where $\operatorname{Tr}$ is the Sobolev trace. Moreover, $\psi\in L^2(\Om)$, and	$f\in L^1(0,T;L^2(\Om))$.
\end{assumption}

At both scales we use the same sign convention for the elliptic operator:
\begin{equation*}
	\Ly w=-\partial_{y_i}(a_{ij}\partial_{y_j}w).
\end{equation*}

\begin{lemma}
	\label{lem:1}
	Let $0<\beta<1$ and let $v:[0,T]\to\R$ be absolutely continuous. For almost every $t\in(0,T)$,
	\begin{equation*}
		v(t)\CD{\beta}v(t)\ge\frac{1}{2}\CD{\beta}\!\left(v^2\right)(t).
	\end{equation*}
\end{lemma}

\begin{lemma}
	\label{lem:2}
	Let $y\ge0$ be absolutely continuous and suppose that
	\[
	\CD{\beta}y(t)\le\lambda y(t)+q(t),\qquad y(0)=y_0,
	\]
	where $\lambda\ge0$ and $q\ge0$ is integrable. Then
	\begin{equation*}
		y(t)\le y_0E_\beta(\lambda t^\beta)
		+\int_0^t(t-s)^{\beta-1}
		E_{\beta,\beta}\!\left(\lambda(t-s)^\beta\right)q(s)\dd s.
	\end{equation*}
	Here $E_{\beta,\nu}(z)=\sum_{m=0}^\infty z^m/\Gamma(\beta m+\nu)$.
	Equivalently, if $a\ge0$ is nondecreasing and
	$y(t)\le a(t)+b(I^\beta y)(t)$, then
	$y(t)\le a(t)E_\beta(bt^\beta)$.
\end{lemma}

Proofs of Lemmas~\ref{lem:1} and \ref{lem:2} can be found in \cite{Alikhanov2010} and \cite{Diethelm2010}, respectively.

\begin{theorem}
\label{thm:3}
Under Assumption~\ref{ass:1}, problem \eqref{eq:1} has a unique weak solution
\[
u^\eps\in C([0,T];H^1(\Om)),\qquad
u_t^\eps\in C([0,T];L^2(\Om)),\qquad
\operatorname{Tr}u^\eps(t)=\widehat u.
\]
Here $\CD{\beta}u_t^\eps:=\frac{\dd}{\dd t}I^{1-\beta}(u_t^\eps-\psi)$
is an $L^2(\Om)$-valued distribution. Multiplication by $\mu^\eps\in L^\infty(\Om)$
is bounded on $L^2(\Om)$, and the weak equation holds in
$\mathcal D'(0,T;H^{-1}(\Om))$. Moreover,
\begin{equation}\label{eq2}
\begin{aligned}
\sup_{0\le t\le T}\bigl(
\norm[L^2(\Om)]{u_t^\eps(t)}^2+\norm[H^1(\Om)]{u^\eps(t)}^2\bigr)
\le C_T\bigl(&\norm[H^1(\Om)]{g}^2+\norm[L^2(\Om)]{\psi}^2\\
&+\norm[L^1(0,T;L^2(\Om))]{f}^2
+\norm[H^{1/2}(\partial\Om)]{\widehat u}^2\bigr),
\end{aligned}
\end{equation}
where $C_T$ is independent of $\eps$.
\end{theorem}
\begin{proof}
Set $V=H_0^1(\Om)$ and, for $v\in H^1(\Om)$ and $\chi\in V$, define
$a^\eps(v,\chi)=\int_\Om A^\eps\nabla v\cdot\nabla\chi\,\dd\bx$.
Let $B^\eps\in H^1(\Om)$ be the  weak harmonic extension
\[
\operatorname{Tr}B^\eps=\widehat u,\qquad
a^\eps(B^\eps,\chi)=0\quad\text{for every }\chi\in V.
\]
A bounded trace lifting and the Lax--Milgram theorem give existence, uniqueness,
and $\norm[H^1(\Om)]{B^\eps}\le C\norm[H^{1/2}(\partial\Om)]{\widehat u}$,
uniformly in $\eps$. Set $w^\eps=u^\eps-B^\eps$. Since
$B_t^\eps=B_{tt}^\eps=\CD{\alpha}B^\eps=0$, the problem for $w^\eps$ is
\[
\begin{cases}
\rho^\eps w_{tt}^\eps+\mu^\eps\CD{\beta}w_t^\eps
-\nabla\!\cdot(A^\eps\nabla w^\eps)=f,\\
\operatorname{Tr}w^\eps=0,\qquad
w^\eps(0)=g-B^\eps\in V,\qquad w_t^\eps(0)=\psi.
\end{cases}
\]

For the resulting homogeneous problem, the proof follows the Volterra reduction and the Galerkin-energy framework in \cite[Theorems~2.4, 2.7, and~2.8]{BakerBanjai2022}, adapted to bounded, spatially heterogeneous coefficients and the weaker regularity assumptions $g - B^\varepsilon \in H_0^1(\Omega)$, $\psi \in L^2(\Omega)$, and $f \in L^1(0, T; L^2(\Omega))$. In this setting, the equation is interpreted distributionally in time, and the Galerkin limit is taken in its time-integrated form. Strong continuity in the energy space follows from approximation by more regular data and the corresponding energy estimate for differences of solutions.   

Finally, $u^\varepsilon = w^\varepsilon + B^\varepsilon$ has the prescribed trace and initial data. The lifting bound combined with the energy estimate of the homogeneous problem directly yields the stated global estimate \eqref{eq2} and continuous dependence on $g, \psi, f$, and $\widehat{u}$.
\end{proof}

\subsection{High-order multi-scale asymptotic analysis}
\label{subsec:2.2}
The chain rule gives the following transformation of the spatial derivatives:
\begin{equation}
	\partial_{x_i}\longmapsto
	\partial_{x_i}+\eps^{-1}\partial_{y_i},
	\qquad i=1,\ldots,d.
	\label{eq:3}
\end{equation}
We assume that $u^\varepsilon(\mathbf{x},t)$ admits the asymptotic expansion 
\begin{equation}
	u^\eps(\bx,t)=u_0(\bx,t,\by)+\eps u_1(\bx,t,\by)
	+\eps^2u_2(\bx,t,\by)+O(\eps^3).
\end{equation}
The material coefficients satisfy
\begin{equation}
	\left\{
	\begin{aligned}	
		a_{ij}^\eps(\bx)&=a_{ij}(\by),\\
		\rho^\eps(\bx)&=\rho(\by),\\
		\mu^\eps(\bx)&=\mu(\by).
	\end{aligned}
	\right.
	\label{eq:5}
\end{equation}
All profiles are $Y$-periodic in $\by$. Substituting \eqref{eq:3}--\eqref{eq:5} into \eqref{eq:1} and collecting terms with equal powers of $\varepsilon$ gives
\begin{align}
	O(\eps^{-2}):&\quad -\partial_{y_i}\!\left(a_{ij}\partial_{y_j}u_0\right)=0, \label{eq:6} \\
	O(\eps^{-1}):&\quad -\partial_{x_i}\!\left(a_{ij}\partial_{y_j}u_0\right)
	-\partial_{y_i}\!\left[a_{ij}\left(\partial_{x_j}u_0+\partial_{y_j}u_1\right)\right]=0, \label{eq:7} \\
	O(\eps^{0}):&\quad \rho u_{0,tt}+\mu\CD{\alpha}u_0
	-\partial_{x_i}\!\left[a_{ij}\left(\partial_{x_j}u_0+\partial_{y_j}u_1\right)\right]\\&\quad 
	-\partial_{y_i}\!\left[a_{ij}\left(\partial_{x_j}u_1+\partial_{y_j}u_2\right)\right]=f. \label{eq:9}
\end{align}
The \(O(\varepsilon^{-2})\)-order equation \eqref{eq:6} implies that \(u_0\) is independent of the microscopic variable \(\mathbf y\). Therefore,
\begin{equation}\label{eq:10}
	u_0(\bx,t,\by)=u_0(\bx,t).
\end{equation}
Taking \eqref{eq:10} into account, all terms involving $\partial_{y_j}u_0$ in the $O(\eps^{-1})$-order equation \eqref{eq:7} vanish. Consequently, \eqref{eq:7} reduces to
\begin{equation*}
	-\partial_{y_i}\left[	
	a_{ij}\left(	
	\partial_{x_j}u_0+\partial_{y_j}u_1	
	\right)\right]=0.	
\end{equation*}
Separation of variables gives the first-order correction in the form
\begin{equation}\label{eq:11}
	u_1(\bx,t,\by)	
	=N^k(\by)\partial_{x_k}u_0(\bx,t),	
\end{equation}
where, for each $k=1,\ldots,d$, the first-order cell function $N^k$ satisfies the auxiliary problem on the PUC $Y$:
\begin{equation}
	\left\{
	\begin{aligned}
		&\Ly N^k=\partial_{y_i}a_{ik} \quad \text{in }Y,\\
		&N^k\text{ is 1-periodic in }\by,\qquad
		\avg{N^k}=0.
	\end{aligned}
	\right.
	\label{eq:12}
\end{equation}
\begin{remark}
	Under Assumption~\ref{ass:1} \textup{(A1)--(A2)}, the bilinear form associated with the first-order auxiliary cell problem \eqref{eq:12} is bounded and coercive on $H_{\#,0}^1(Y)$. Moreover, its right-hand side defines a bounded linear functional on this space. Therefore, the Lax--Milgram theorem ensures that, for each $k=1,\ldots,d$, the problem \eqref{eq:12} admits a unique weak solution $N^k\in H_{\#,0}^1(Y)$.
\end{remark}
Substituting \eqref{eq:11} into the $O(\eps^0)$-order equation \eqref{eq:9}, integrating over $Y$, and applying the divergence theorem yields the macroscopic homogenized problem
\begin{equation}
	\left\{
	\begin{aligned}
		&\rhoh u_{0,tt}
		+\muh\CD{\alpha}u_0
		-\partial_{x_k}
		\left(
		\widehat a_{k\ell}
		\partial_{x_\ell}u_0
		\right)
		=f
		&&\text{in }\Om\times(0,T],\\
		&u_0=\widehat u(\bx)
		&&\text{on }\partial\Om\times[0,T],\\
		&u_0(\bx,0)=g(\bx)
		&&\text{in }\Om,\\
		&u_{0,t}(\bx,0)=\psi(\bx)
		&&\text{in }\Om.
	\end{aligned}
	\right.
	\label{eq:13}
\end{equation}

Here, the homogenized material parameters in \eqref{eq:13} are obtained by integration over the PUC $Y$:
\begin{equation*}
	\rhoh=\avg{\rho},
	\qquad
	\muh=\avg{\mu},
	\qquad
	\widehat a_{k\ell}
	=
	\avg{
		a_{k\ell}
		+a_{kj}\partial_{y_j}N^\ell
	}.
\end{equation*}
We next derive the second-order correction term $u_2$. Substituting \eqref{eq:10} and \eqref{eq:11} into \eqref{eq:9} gives
\begin{align}
	\Ly u_2={}&(\rhoh-\rho)u_{0,tt}+(\muh-\mu)\CD{\alpha}u_0\notag\\
	&+\left[
	a_{k\ell}+a_{kj}\partial_{y_j}N^\ell
	+\partial_{y_i}(a_{ik}N^\ell)-\widehat a_{k\ell}
	\right]\partial_{x_kx_\ell}u_0.
	\label{eq:14}
\end{align}
Equation~\eqref{eq:14} gives the second-order correction term in the form
\begin{equation}
	u_2(\bx,\by,t)=F(\by)u_{0,tt}(\bx,t)+G(\by)\,{}^C D_t^{\alpha} u_0(\bx,t)
	+N^{k\ell}(\by)\partial_{x_kx_\ell}u_0(\bx,t).
	\label{eq:15}
\end{equation}
Applying $\Ly$ to \eqref{eq:15} and comparing coefficients in \eqref{eq:14} defines the second-order auxiliary cell problems
\begin{subequations}
	\begin{align}
		\left\{
		\begin{aligned}
			&\Ly F=\rhoh-\rho \quad \text{in }Y,\\
			&F\text{ is 1-periodic in }\by,\qquad \avg{F}=0,
		\end{aligned}
		\right.
		\label{eq:16a}\\
		\left\{
		\begin{aligned}
			&\Ly G=\muh-\mu\quad \text{in }Y,\\
			&G\text{ is 1-periodic in }\by,\qquad \avg{G}=0,
		\end{aligned}
		\right.
\\
		\left\{
		\begin{aligned}
			&\Ly N^{k\ell}
			=a_{k\ell}+a_{kj}\partial_{y_j}N^\ell
			+\partial_{y_i}(a_{ik}N^\ell)-\widehat a_{k\ell}
			\quad \text{in }Y,\\
			&N^{k\ell}\text{ is 1-periodic in }\by,\qquad \avg{N^{k\ell}}=0,
			\qquad 1\le k,\ell\le d.
		\end{aligned}
		\right.
		\label{eq:16c}
	\end{align}
\end{subequations}

\begin{remark}
	The source terms in \eqref{eq:16a}--\eqref{eq:16c}, namely \(\hat{\rho}-\rho\), \(\hat{\mu}-\mu\), and \(a_{k\ell}+a_{kj}\partial_{y_j}N^\ell+\partial_{y_i}(a_{ik}N^\ell)-\hat{a}_{k\ell}\), have zero average over the PUC \(Y\) by the definitions of the homogenized material parameters. Under Assumption~\ref{ass:1} \textup{(A1)--(A2)}, the Lax--Milgram theorem guarantees existence and uniqueness of solutions to \eqref{eq:16a}--\eqref{eq:16c} for any fixed macroscopic field \(u_0\).
\end{remark}
The homogenized, low-order multi-scale (LOMS), and HOMS asymptotic solutions to the fractionally damped wave problem \eqref{eq:1} are therefore

\begin{align}
	u^{(0)}(\bx,t)
	=&u_0(\bx,t),\notag
	\\
	u^{(1\varepsilon)}(\bx,t)
	=&u_0(\bx,t)+\eps N^k(\by)\partial_{x_k}u_0(\bx,t),\notag
	\\
	u^{(2\varepsilon)}(\bx,t)
	=&u_0(\bx,t)+\eps N^k(\by)\partial_{x_k}u_0(\bx,t)
	+\eps^2[F(\by)u_{0,tt}(\bx,t)\notag\\&+G(\by)\CD{\alpha}u_0(\bx,t)+N^{k\ell}(\by)\partial_{x_kx_\ell}u_0(\bx,t)].
	\label{eq:17}
\end{align}

The cell functions in \eqref{eq:17} recover microscopic variations absent from the homogenized solution: $N^k$ restores stiffness-induced gradient and flux fluctuations; $F$ and $G$ correct local inertia and history-dependent damping associated with density and damping-coefficient fluctuations, respectively; and $N^{k\ell}$ corrects stiffness-flux variations driven by the macroscopic curvature. The LOMS solution includes the first-order $N^k$ contribution, while the HOMS solution additionally incorporates the inertia, fractional-damping, and curvature corrections. As shown in Section~\ref{sec:3}, these second-order terms cancel the leading $O(1)$ residual source remaining in the LOMS solution, yielding a formal $O(\eps)$ local pointwise strong residual under the stated regularity conditions. This improved local balance explains the benefit of the HOMS solution.

\section{Error analysis of HOMS asymptotic solutions}
\label{sec:3}
\subsection{Error analysis of HOMS solutions in the pointwise sense}
\begin{assumption}
	\label{ass:2}
	The homogenized solution is assumed to be sufficiently regular. One explicit sufficient condition is
	\begin{equation}
		u_0\in W^{4,2}(0,T;W^{4,\infty}(\Om)),
		\qquad u_{0,tt}(0)=u_{0,ttt}(0)=0.
        \label{eq:18}
    \end{equation}
\end{assumption}

We first define the residual functions, namely the differences between the exact solution and the LOMS and HOMS asymptotic solutions:
\begin{equation*}
	\begin{cases} u_{\Delta}^{(1\varepsilon)}(\boldsymbol{x}, t) = u^{\varepsilon}(\boldsymbol{x}, t) - u^{(1\varepsilon)}(\boldsymbol{x}, t),\\ u_{\Delta}^{(2\varepsilon)}(\boldsymbol{x}, t) = u^{\varepsilon}(\boldsymbol{x}, t) - u^{(2\varepsilon)}(\boldsymbol{x}, t). \end{cases} 
\end{equation*}
Substituting the definition of $u_{\Delta}^{(1\varepsilon)}$ into \eqref{eq:1} gives the residual equations for the LOMS solution
\begin{equation}
	\left\{\begin{aligned}
		&\rho^\eps u_{\Delta,tt}^{(1\varepsilon)}
		+\mu^\eps\CD{\alpha}u_{\Delta}^{(1\varepsilon)}\\
		&\quad-\partial_{x_i}\!\left(a_{ij}^\eps\partial_{x_j}u_{\Delta}^{(1\varepsilon)}\right)
		=F_0^\eps(\bx,\by,t)+\eps F_1^\eps(\bx,\by,t),
		&&\text{in }\Om\times(0,T],\\[2pt]
		&u_{\Delta}^{(1\varepsilon)}(\bx,t)
		=-\eps N^k\!\left(\by\right)
		\partial_{x_k}u_0(\bx,t),
		&&\text{on }\partial\Om\times[0,T],\\[2pt]
		&u_{\Delta}^{(1\varepsilon)}(\bx,0)
		=-\eps N^k\!\left(\by\right)
		\partial_{x_k}g(\bx),
		&&\text{in }\Om,\\[2pt]
		&u_{\Delta,t}^{(1\varepsilon)}(\bx,0)
		=-\eps N^k\!\left(\by\right)
		\partial_{x_k}\psi(\bx),
		&&\text{in }\Om,
	\end{aligned}\right.
\end{equation}
where
\begin{align*}
	F_0^\eps={}&
	(\rhoh-\rho)u_{0,tt}
	+(\muh-\mu)\CD{\alpha}u_0+
	\left[
	a_{k\ell}
	+a_{kj}\partial_{y_j}N^\ell
	+\partial_{y_i}(a_{ik}N^\ell)
	-\widehat a_{k\ell}
	\right]\partial_{x_kx_\ell}u_0,
	\\
	F_1^\eps={}&
	-\rho N^k \partial_{x_k}u_{0,tt}
	-\mu N^k\CD{\alpha}
	\left(\partial_{x_k}u_0\right)
	+a_{ij}N^k \partial_{x_ix_jx_k}u_0.
\end{align*}
Similarly, substituting the residual function $u_{\Delta}^{(2\varepsilon)}(\boldsymbol{x}, t)$ into \eqref{eq:1} gives the residual equations for the HOMS solution
\begin{equation}
	\left\{\begin{aligned}
		&\rho^\eps u_{\Delta,tt}^{(2\varepsilon)}
		+\mu^\eps\CD{\alpha}u_{\Delta}^{(2\varepsilon)}
		-\partial_{x_i}\!\left(a_{ij}^\eps\partial_{x_j}u_{\Delta}^{(2\varepsilon)}\right)
		=\eps F_2^\eps(\bx,\by,t),
		&&\text{in }\Om\times(0,T],\\[2pt]
		&u_{\Delta}^{(2\varepsilon)}(\bx,t)
		=-\eps N^k(\by)\partial_{x_k}u_0(\bx,t)
		-\eps^2\Bigl[F(\by)u_{0,tt}(\bx,t)\\
		&\qquad+G(\by)\CD{\alpha}u_0(\bx,t)
		+N^{k\ell}(\by)\partial_{x_kx_\ell}u_0(\bx,t)\Bigr],
		&&\text{on }\partial\Om\times[0,T],\\[2pt]
		&u_{\Delta}^{(2\varepsilon)}(\bx,0)
		=-\eps N^k(\by)\partial_{x_k}g(\bx)
		-\eps^2N^{k\ell}(\by)\partial_{x_kx_\ell}g(\bx),
		&&\text{in }\Om,\\[2pt]
		&u_{\Delta,t}^{(2\varepsilon)}(\bx,0)
		=-\eps N^k(\by)\partial_{x_k}\psi(\bx)
		-\eps^2N^{k\ell}(\by)\partial_{x_kx_\ell}\psi(\bx),
		&&\text{in }\Om.
	\end{aligned}\right.
\end{equation}
The residual source is
\begin{align*}
	F_2^\eps={}&F_1^\eps+a_{ij}\partial_{y_j}\partial_{x_i}u_2
	+\partial_{y_i}(a_{ij}\partial_{x_j}u_2)\\
	&+\eps\left[a_{ij}\partial_{x_i x_j}u_2
	-\rho\partial_{tt}u_2-\mu\CD{\alpha}u_2\right].
\end{align*}

The above pointwise error results clearly indicate that the LOMS solution fails to preserve local physical balance because the source term $F_0^\varepsilon$ in its residual equation, which carries no prefactor of $\varepsilon$, does not vanish as the microstructural parameter $\varepsilon$ tends to zero, resulting in a formal $O(1)$ strong residual source. By incorporating second-order correction terms, the HOMS solution achieves an $O(\varepsilon)$ pointwise error and cancels the leading local balance defect remaining in the LOMS solution, thereby substantially improving and ensuring local physical balance. This provides the principal motivation for constructing an HOMS solution that enables highly accurate computations for composite materials.

\begin{remark}
	\label{rem:3}
	Let $\gamma=2-\alpha$. For $t>0$, differentiating
	\(	\CD{\alpha}u_0=I^\gamma u_{0,tt}\)
	twice gives
	\begin{equation*}
		\partial_{tt}\left(\CD{\alpha}u_0\right)
		=
		I^\gamma u_{0,tttt}
		+\frac{u_{0,ttt}(0)t^{1-\alpha}}{\Gamma(2-\alpha)}
		+\frac{(1-\alpha)u_{0,tt}(0)t^{-\alpha}}{\Gamma(2-\alpha)}.
	\end{equation*}
	The factor $t^{-\alpha}$ does not belong to
	$L^1(0,\delta)$ for any $\delta>0$ when $\alpha>1$.
	Hence, if $u_{0,tt}(0)\ne0$, then
	\(	\CD{\alpha}u_0
	\notin W^{2,1}(0,T;H^2(\Om)),
	\)
	and the iterated Caputo derivative
	\(
	\CD{\alpha}\left(\CD{\alpha}u_0\right)
	=
	I^\gamma
	\partial_{tt}\left(\CD{\alpha}u_0\right)\)	is not defined in the standard sense. If $u_{0,tt}(0)=0$ but $u_{0,ttt}(0)\ne0$, the remaining endpoint term $t^{1-\alpha}$ belongs to $L^1(0,\delta)$ for every $1<\alpha<2$, and belongs to $L^2(0,\delta)$ exactly when $\alpha<3/2$. Thus, the two endpoint conditions in \eqref{eq:18} provide a simple sufficient condition for an unweighted $L^2$ residual source for the HOMS solution over the full range $1<\alpha<2$. When $\alpha<3/2$, the condition $u_{0,ttt}(0)=0$ may be relaxed.
\end{remark}

\subsection{Error analysis of HOMS solution in the integral sense}
\label{subsec:3.2}
To ensure the computational reliability of global wave simulations, we provide a global error estimate in the integral sense for the high-order multi-scale solutions. To this end, we first list the additional assumptions required for the subsequent proof.
\begin{assumption}
	\label{ass:3}
	\par\smallskip\noindent\textup{(B1)}\quad
	The domain consists of complete periodicity cells, namely
	\[
	\overline\Om=\bigcup_{z\in I_\eps}\eps(z+\overline Y),\qquad
	I_\eps=\{z\in\mathbb Z^d:\eps(z+\overline Y)\subset\overline\Om\}.
	\]
	Write $E_z=\eps(z+Y)$ for a translated cell and $\partial E_z$ for its boundary.
	\par\smallskip\noindent\textup{(B2)}\quad
	Let $a_{ij}(\by)=a_{ii}(\by)\delta_{ij}$, where $\delta_{ij}$ is the Kronecker symbol. The material coefficients $a_{ii}$, $\rho$, and $\mu$ are symmetric with respect to every middle hyperplane of $Y=[0,1]^d$.
	\par\smallskip\noindent\textup{(B3)}\quad
	All auxiliary cell functions are subject to homogeneous Dirichlet boundary conditions instead of periodic boundary conditions \cite{Cao2006,CaoCuiZhu2002,WangCaoWong2015}.
\end{assumption}

\begin{lemma}
	\label{lem:4}
	Under \textup{(A1)}--\textup{(A2)} and \textup{(B2)}--\textup{(B3)}, the normal flux $n_i a_{ij}(\by)\partial_{y_j}\chi$ of every first- and second-order cell function $\chi$ is continuous across opposite faces of the PUC $Y$. This follows from the corresponding cell-boundary argument in \cite{Cao2006,CaoCuiZhu2002,WangCaoWong2015}.
\end{lemma}

Set $e^\eps=u^\eps-u^{(2\eps)}$. Write the initial errors as
\begin{equation}
	\begin{aligned}
		e^\eps(\bx,0)&=:\eps\widetilde\eta_2(\bx),&
		\widetilde\eta_2&=-N^{k,\eps}\partial_{x_k}g-\eps N^{k\ell,\eps}\partial_{x_kx_\ell}g,\\
		e_t^\eps(\bx,0)&=:\eps\widetilde\zeta_2(\bx),&
		\widetilde\zeta_2&=-N^{k,\eps}\partial_{x_k}\psi-\eps N^{k\ell,\eps}\partial_{x_kx_\ell}\psi.
	\end{aligned}
\end{equation}
With $P^\eps=P_1^\eps+\eps P_2^\eps$, where $P_1^\eps$, $P_2^\eps$, and $\bm H^\eps$ are specified in the proof below, the error $e^\eps$ satisfies the following initial-boundary value problem, with the differential equation understood in the spatial weak sense:
\begin{equation}
	\left\{\begin{aligned}
		&\rho^\eps e_{tt}^\eps+\mu^\eps\CD{\beta}e_t^\eps
		-\nabla\cdot(A^\eps\nabla e^\eps)
		=\eps P^\eps+\eps^2\nabla\cdot\bm H^\eps,
		&&\text{in }\Om\times(0,T],\\
		&e^\eps=0,&&\text{on }\partial\Om\times[0,T],\\
		&e^\eps(\bx,0)=\eps\widetilde\eta_2(\bx),&&\text{in }\Om,\\
		&e_t^\eps(\bx,0)=\eps\widetilde\zeta_2(\bx),&&\text{in }\Om.
	\end{aligned}\right.
	\label{eq:22}
\end{equation}
Define
\begin{equation}
	\mathcal E^\eps(t)=\tfrac12\norm[L^2(\Omega)]{\sqrt{\rho^\eps}e_t^\eps(t)}^2
	+\tfrac12(A^\eps\nabla e^\eps(t),\nabla e^\eps(t)),
	\qquad M^\eps(t)=\norm[L^2(\Omega)]{\sqrt{\mu^\eps}e_t^\eps(t)}^2.
\end{equation}

\begin{theorem}
	\label{thm:5}
	Let Assumptions~\ref{ass:1}, \ref{ass:2}, and \ref{ass:3} hold.
	Then
	\begin{align}
		&\sup_{0\le t\le T}\left(
		\norm[L^2(\Omega)]{e_t^\eps(t)}^2+\norm[H_0^1(\Omega)]{e^\eps(t)}^2\right)
		+\sup_{0\le t\le T}(I^{1-\beta}M^\eps)(t)
		\le C_T\eps^2.
		\label{eq:24}
	\end{align}
	Consequently,
	\begin{equation}
		\sup_{0\le t\le T}\left(
		\norm[L^2(\Omega)]{u_t^\eps-u_t^{(2\eps)}}
		+\norm[H_0^1(\Omega)]{u^\eps-u^{(2\eps)}}\right)\le C_T\eps.
		\label{eq:25}
	\end{equation}
	The constants are independent of $\eps$.
\end{theorem}
\begin{proof}
	Substituting \eqref{eq:15} into the residual source and collecting terms gives
	$\eps F_2^\eps=\eps P^\eps+\eps^2\nabla_{\bx}\cdot\bm H^\eps$,
	where $P^\eps=P_1^\eps+\eps P_2^\eps$ and
	\begin{align*}
		P_1={}&
		-\rho N^k \partial_{x_k}u_{0,tt}
		-\mu N^k\partial_{x_k}
		\left(\CD{\alpha}u_0\right)
		+a_{ij}N^k \partial_{x_ix_jx_k}u_0
		\notag\\
		&+a_{ij}\left[
		(\partial_{y_j}F)\partial_{x_i}u_{0,tt}
		+(\partial_{y_j}G)
		\partial_{x_i}\left(\CD{\alpha}u_0\right)
		+(\partial_{y_j}N^{k\ell})\partial_{x_ix_kx_\ell}u_0
		\right],
		\\
		P_2={}&
		-\rho\left[
		Fu_{0,tttt}
		+G\partial_{tt}\left(\CD{\alpha}u_0\right)
		+N^{k\ell}\partial_{x_kx_\ell}u_{0,tt}
		\right]
		\notag\\
		&-\mu\left[
		F\CD{\alpha}(u_{0,tt})
		+G\CD{\alpha}\left(\CD{\alpha}u_0\right)
		+N^{k\ell}\CD{\alpha}(\partial_{x_kx_\ell}u_0)
		\right],
		\\
		H_i={}&
		a_{ij}\left[
		F\partial_{x_j}u_{0,tt}
		+G\partial_{x_j}\left(\CD{\alpha}u_0\right)
		+N^{k\ell}\partial_{x_jx_kx_\ell}u_0
		\right].
	\end{align*}
	The regularity in \eqref{eq:18}, the cell-function $H^1(Y)$ bounds, and scaling over the complete cells give
	\begin{equation*}
		\begin{aligned}
			&\norm[L^2(0,T;L^2(\Om))]{P_1^\eps}
			+\norm[L^2(0,T;L^2(\Om))]{P_2^\eps}\\
			&\quad+\norm[L^\infty(0,T;L^2(\Om))]{\bm H^\eps}
			+\norm[L^2(0,T;L^2(\Om))]{\bm H_t^\eps}\le C,
		\end{aligned}
	\end{equation*}
	where $C$ is independent of $\eps$. The residual identity is understood globally in $L^2(0,T;H^{-1}(\Om))$ when performing the energy calculation.

	The available regularity, \(e_t^\varepsilon \in L^\infty(0,T; L^2(\Omega))\), does not justify using \(e_t^\varepsilon\) directly as a test function. We justify the following calculation by the density argument in \cite{CioranescuDonato1999,DongCao2009}  and omit the details for brevity. We test the equation in \eqref{eq:22} with \(e_t^\varepsilon\), integrate over \(\Omega\), and apply Green's formula with the boundary condition to obtain
	\begin{align*}
		(\rho^\varepsilon e_{tt}^\varepsilon, e_t^\varepsilon) &= \frac{1}{2} \frac{\mathrm{d}}{\mathrm{d}t} \int_\Omega \rho^\varepsilon |e_t^\varepsilon|^2 \, \mathbf{dx}, \\
		(A^\varepsilon \nabla e^\varepsilon, \nabla e_t^\varepsilon) &= \frac{1}{2} \frac{\mathrm{d}}{\mathrm{d}t} \int_\Omega A^\varepsilon \nabla e^\varepsilon \cdot \nabla e^\varepsilon \, \mathbf{dx}.
	\end{align*}
	Applying Lemma~\ref{lem:1} pointwise in $\boldsymbol{x}$, multiplying by $\mu^\varepsilon$, and integrating over $\Omega$ gives
	\[
	(\mu^\varepsilon {}^C D_t^\beta e_t^\varepsilon, e_t^\varepsilon) \geq \frac{1}{2} {}^C D_t^\beta M^\varepsilon(t).
	\]
	Consequently,
	\[
	\frac{\dd}{\dd t}\mathcal E^\eps+\tfrac12\CD{\beta}M^\eps
	\le\eps(P^\eps,e_t^\eps)-\eps^2(\bm H^\eps,\nabla e_t^\eps).
	\]
	Using
	\[
	\int_0^t\CD{\beta}M^\eps(s)\dd s
	=(I^{1-\beta}M^\eps)(t)
	-\frac{t^{1-\beta}}{\Gamma(2-\beta)}M^\eps(0),
	\]
	integrating in time, and integrating the flux term by parts in time yield
	\begin{align}
		\mathcal E^\eps(t)+\tfrac12(I^{1-\beta}M^\eps)(t)
		\le{}&\mathcal E^\eps(0)
		+\frac{t^{1-\beta}}{2\Gamma(2-\beta)}M^\eps(0)
		+\eps\int_0^t(P^\eps,e_t^\eps)\dd s\notag\\
		&-\eps^2(\bm H^\eps(t),\nabla e^\eps(t))
		+\eps^2(\bm H^\eps(0),\nabla e^\eps(0))\notag\\
		&+\eps^2\int_0^t(\bm H_t^\eps,\nabla e^\eps)\dd s.
		\label{eq:26}
	\end{align}
	Next, set $\mathcal Y^\eps=\mathcal E^\eps+\tfrac12 I^{1-\beta}M^\eps$.
	Substituting the initial conditions into the energy terms gives
	\begin{equation*}
		\begin{aligned}
			\mathcal E^\eps(0)+M^\eps(0)
			={}&\tfrac12\norm[L^2(\Omega)]{\sqrt{\rho^\eps}\,\eps\widetilde\zeta_2}^2
			+\tfrac12\bigl(A^\eps\nabla(\eps\widetilde\eta_2),
			\nabla(\eps\widetilde\eta_2)\bigr)\\
			&+\norm[L^2(\Omega)]{\sqrt{\mu^\eps}\,\eps\widetilde\zeta_2}^2
			\le C\eps^2.
		\end{aligned}
	\end{equation*}
	The Cauchy--Schwarz and Young inequalities, together with the residual-component estimate above, then give
	\begin{align*}
		\eps\int_0^t|(P^\eps,e_t^\eps)|\dd s
		&\le C\eps^2+C\int_0^t\mathcal Y^\eps(s)\dd s,\\
		\eps^2|(\bm H^\eps(t),\nabla e^\eps(t))|
		&\le\tfrac14\mathcal E^\eps(t)+C\eps^4,\\
		\eps^2|(\bm H^\eps(0),\nabla e^\eps(0))|&\le C\eps^3,\\
		\eps^2\int_0^t|(\bm H_t^\eps,\nabla e^\eps)|\dd s
		&\le C\eps^4+C\int_0^t\mathcal Y^\eps(s)\dd s.
	\end{align*}
	Absorbing the energy term in \eqref{eq:26}, using $\eps\le1$, and applying the integral Gronwall inequality give
	$\mathcal Y^\eps(t)\le C_T\eps^2$.
	Consequently, ellipticity, the Poincar\'e inequality, and the nonnegativity of the memory term prove
	\eqref{eq:24}. Taking square roots gives \eqref{eq:25}.
\end{proof}

\section{Multi-scale numerical algorithm}
\label{sec:4}

\subsection{Computation for microscopic cell functions and homogenized material parameters}
\begin{enumerate}

	\item[(1)] On the PUC $Y=(0,1)^d$, generate shape-regular, interface-fitted triangular or tetrahedral meshes
	\(	\mathcal T_{h_Y}=\{K\},\) and set \(
	h_Y=\max_{K\in\mathcal T_{h_Y}}\operatorname{diam}(K).\)
	Define the conforming finite element space
	\[	S_{h_Y}
	=
	\left\{
	v_h\in H_0^1(Y):
	v_h|_K\in\mathbb P_{p_Y}(K),
	\quad K\in\mathcal T_{h_Y}
	\right\}.\]

	\item[(2)] 	All numerical auxiliary cell functions satisfy homogeneous Dirichlet boundary conditions \cite{Cao2006,CaoCuiZhu2002,WangCaoWong2015}. For $v,w\in H_0^1(Y)$, introduce the cell bilinear form
	\(	\mathfrak a_Y(v,w)	=	\int_Y	A\nabla_{\by}v\cdot\nabla_{\by}w 	\dd\by.\)
	For each $1\le k\le d$, find	$N_{h_Y}^k\in S_{h_Y}$ such that
	\begin{equation}
		\mathfrak a_Y(N_{h_Y}^k,\phi_h)
		=
		-\int_Y
		a_{ik}\partial_{y_i}\phi_h
		\dd\by,
		\qquad
		\forall\phi_h\in S_{h_Y}.
		\label{eq:27}
	\end{equation}

	Using the computed first-order cell functions, evaluate the	homogenized material parameters as
	\begin{equation}
		\rhoh=\avg{\rho},
		\qquad
		\muh=\avg{\mu},
		\qquad
		\widehat a_{k\ell,h_Y}
		=
		\avg{
			a_{k\ell}
			+a_{kj}\partial_{y_j}N_{h_Y}^{\ell}
		}.
		\label{eq:28}
	\end{equation}
	\item[(3)] 	Using the same mesh and finite element space, find
	$F_{h_Y}$, $G_{h_Y}$, and
	$N_{h_Y}^{k\ell}\in S_{h_Y}$ such that
	\begin{subequations}
		\label{eq:29}
		\begin{align}
			\mathfrak a_Y(F_{h_Y},\phi_h)
			=&
			\int_Y(\rhoh-\rho)\phi_h\dd\by,
			\\
			\mathfrak a_Y(G_{h_Y},\phi_h)
			=&
			\int_Y(\muh-\mu)\phi_h\dd\by,
			\\
			\mathfrak a_Y(N_{h_Y}^{k\ell},\phi_h)
			=&
			\int_Y
			\left(
			a_{k\ell}
			+a_{kj}\partial_{y_j}N_{h_Y}^{\ell}
			-\widehat a_{k\ell,h_Y}
			\right)
			\phi_h\dd\by\notag
			\\&-\int_Y
			a_{ik}N_{h_Y}^{\ell}
			\partial_{y_i}\phi_h
			\dd\by,
		\end{align}
	\end{subequations}
	for every $\phi_h\in S_{h_Y}$ and
	$1\le k,\ell\le d$.	
\end{enumerate}

\subsection{Computation for macroscopic homogenized equation by a fully discrete \texorpdfstring{$L2$-$1_\sigma$}{L2-1sigma}--Newmark finite element scheme}
\label{subsec:fully-discrete-scheme}

Let $p\ge1$, and let $\mathcal T_h$ be a conforming, shape-regular,
quasi-uniform simplicial mesh of $\Om$. Define the conforming finite element test space by

\[
V_h=\{v_h\in C^0(\overline\Om):
v_h|_K\in\mathbb P_p(K),\ K\in\mathcal T_h,
\ v_h|_{\partial\Om}=0\}\subset H_0^1(\Om).
\]

In the numerical experiments, we take $p=2$. Assemble the constant bilinear forms
	\[
	m(w,z)=\rhoh(w,z),\qquad
	d(w,z)=\muh(w,z),\qquad
	a_h(w,z)=(\Ah_{h_Y}\nabla w,\nabla z),
	\]
	where $\rhoh$ and $\muh$ are the cell averages and $\Ah_{h_Y}$
	is computed from \eqref{eq:27}--\eqref{eq:28}.

For the $L2$-$1_\sigma$ approximation \cite{Alikhanov2015}, set $\beta=\alpha-1\in(0,1)$ and $\sigma=1-\beta/2=(3-\alpha)/2$.
Writing $v=u_{0,t}$, the macroscopic homogenized equation becomes
\begin{equation*}
	\rhoh v_t+\muh\CD{\beta}v-\nabla\cdot(\Ah\nabla u_0)=f,
	\qquad u_{0,t}=v.
\end{equation*}
Let $t_n=n\tau$, $0\le n\le N$, $\tau=T/N$, and $t_{n+\sigma}=(n+\sigma)\tau$. For a grid function $\{W^n\}$, define
\begin{equation*}
	(\Ddisc_\tau^\beta W)^{n+\sigma}
	=
	\frac{\tau^{-\beta}}{\Gamma(2-\beta)}
	\sum_{k=0}^{n} C_{n-k}^{(n+1)}(W^{k+1}-W^k).
\end{equation*}
Here, the auxiliary weights are
\begin{align*}
	A_0=&\sigma^{1-\beta},\\
	A_j=&(j+\sigma)^{1-\beta}-(j-1+\sigma)^{1-\beta},
	\qquad j\ge1,\\
	B_j=&
	\frac{(j+\sigma)^{2-\beta}-(j-1+\sigma)^{2-\beta}}{2-\beta}\\&
	-\frac{1}{2}\left[(j+\sigma)^{1-\beta}+(j-1+\sigma)^{1-\beta}\right],
	\quad j\ge1.
\end{align*}
The convolution weights are $C_0^{(1)}=A_0$ and, for $n\ge1$,
\begin{equation}
	\begin{aligned}
		C_0^{(n+1)}&=A_0+B_1,\\
		C_j^{(n+1)}&=A_j+B_{j+1}-B_j,\qquad 1\le j\le n-1,\\
		C_n^{(n+1)}&=A_n-B_n.
	\end{aligned}
	\label{eq:30}
\end{equation}
We next introduce the discrete macroscopic approximations
\[
U_h^n\approx u_0(t_n),\qquad
V_h^n\approx u_{0,t}(t_n),\qquad
A_h^n\approx u_{0,tt}(t_n).
\]
The Newmark updates \cite{Newmark1959} with $\gamma_N=1/2$ and $\beta_N=1/4$ are
\begin{align}
	U_h^{n+1}
	&=U_h^n+\tau V_h^n
	+\frac{\tau^2}{4}\left(A_h^n+A_h^{n+1}\right),
	\label{eq:31}\\
	V_h^{n+1}
	&=V_h^n
	+\frac{\tau}{2}\left(A_h^n+A_h^{n+1}\right).
	\label{eq:32}
\end{align}

Define the shifted values and fractional derivative by
\begin{equation*}
	\begin{aligned}
		Z_h^{n+\sigma}=\sigma Z_h^{n+1}+(1-\sigma)Z_h^n,
		\qquad Z\in\{U,V,A\},\\
		Q_h^{n+\sigma}
		:=(\Ddisc_\tau^\beta V_h)^{n+\sigma}
		\approx \CD{\beta}u_{0,t}(t_{n+\sigma})
		=\CD{\alpha}u_0(t_{n+\sigma}).
	\end{aligned}
\end{equation*}
The fully discrete homogenized problem is to find
a continuous piecewise polynomial displacement $U_h^{n+1}$ of degree at most $p$ on each element, and $V_h^{n+1},A_h^{n+1}\in V_h$, satisfying
\begin{equation}
\left\{\begin{aligned}
&m(A_h^{n+\sigma},\chi_h)+d(Q_h^{n+\sigma},\chi_h)
+a_h(U_h^{n+\sigma},\chi_h)\\
&\qquad=(f^{n+\sigma},\chi_h),\qquad\forall\chi_h\in V_h,\\
&U_h^{n+1}=\widehat u_h\qquad\text{on }\partial\Om.
\end{aligned}\right.
\label{eq:33}
\end{equation}

Let $I_h$ be the degree-$p$ Lagrange interpolation operator on $\mathcal T_h$ and set $\widehat u_h=\operatorname{Tr}(I_hg)$. For the convergence analysis, initialize
\[
U_h^0=R_h^Dg,\qquad V_h^0=P_h\psi,
\]
and determine \(A_h^0\in V_h\) from
\begin{equation}
	m(A_h^0,\chi_h)+a_h(U_h^0,\chi_h)
	=(f(0),\chi_h)
	\quad\forall\chi_h\in V_h.
	\label{eq:34}
\end{equation}
Here $P_h:L^2(\Om)\to V_h$ is the $L^2$ projection. For $w\in H_0^1(\Om)$,
the homogeneous Ritz projection $R_hw\in V_h$ is defined by
\[
a_h(R_hw,\chi_h)=a_h(w,\chi_h),\qquad\chi_h\in V_h.
\]
For $w\in H^1(\Om)$ with $\operatorname{Tr}w=\widehat u$, define the degree-$p$ conforming Ritz projection $R_h^Dw$ by
\[
\begin{cases}
a_h(R_h^Dw,\chi_h)=a_h(w,\chi_h),&\chi_h\in V_h,\\
R_h^Dw=\widehat u_h,&\text{on }\partial\Om.
\end{cases}
\]

\subsection{Computation and reconstruction for the HOMS solution}
At each step, eliminate $V_h^{n+1}$ and $A_h^{n+1}$ from \eqref{eq:33} using \eqref{eq:31}~and~\eqref{eq:32}, solve for $U_h^{n+1}$, and recover the velocity and acceleration. Evaluate the shifted fields and $Q_h^{n+\sigma}$ as defined above.

The homogenized, LOMS, and HOMS reconstructions at $t_{n+\sigma}$ are
	\begin{align*}
		u_h^{[0],n+\sigma}=&U_h^{n+\sigma},\\
		u_h^{[1],n+\sigma}=&U_h^{n+\sigma}
		+\eps N_{h_Y}^k(\by)\,\partial_{x_k}U_h^{n+\sigma},\\
		u_h^{[2],n+\sigma}=&u_h^{[1],n+\sigma}
		+\eps^2\left[F_{h_Y}(\by)A_h^{n+\sigma}\right.\\
		&\left.+G_{h_Y}(\by)Q_h^{n+\sigma}
		+N_{h_Y}^{k\ell}(\by)\,\partial_{x_kx_\ell}U_h^{n+\sigma}\right].
	\end{align*}

\section{Stability and error analysis  for multi-scale algorithm}
\label{sec:5}

For the proposed multi-scale algorithm, we analyze the stability of the fully discrete macroscopic scheme and the errors arising from discretization at the microscopic cell and macroscopic scales. In the numerical analysis, $\Ah$ denotes the macroscopic homogenized material parameter defined by the exact continuous cell problems, and $u_0$ denotes the corresponding continuous macroscopic solution. All constants may depend on $T$ and $\beta$, but are independent of $h_Y$, $h$, and $\tau$.

\subsection{Shifted-energy stability for the fully discrete scheme at macro-scale}

Eliminating the velocity and acceleration from \eqref{eq:33} gives the bilinear form
\[
\frac{4\sigma}{\tau^2}m(w_h,\chi_h)
+\frac{2\tau^{-\beta}C_0^{(n+1)}}{\tau\Gamma(2-\beta)}d(w_h,\chi_h)
+\sigma a_h(w_h,\chi_h).
\]
The tensor $\Ah_{h_Y}$ is symmetric positive definite, with bounds independent of $h_Y$. Together with $C_0^{(n+1)}>0$ and $\rhoh,\muh>0$, this makes the bilinear form coercive.  Hence every time step has a unique solution by the Lax--Milgram theorem, and the Newmark updates give unique velocities and accelerations in $V_h$.
For the ensuing stability and consistency analysis, note that when $\sigma=1-\beta/2$, the weights are positive and decrease with the lag \cite[Lemma~4]{Alikhanov2015}. If $W\in C^3([0,T];L^2(\Omega))$, then
\begin{equation}
	\norm[L^2(\Omega)]{\CD{\beta}W(t_{n+\sigma})
		-(\Ddisc_\tau^\beta W)^{n+\sigma}}
	\le C\tau^{3-\beta}
	\max_{0\le t\le t_{n+1}}\norm[L^2(\Omega)]{W^{(3)}(t)}.
	\label{eq:35}
\end{equation}
This is the $L^2(\Omega)$-valued form of the local consistency estimate in \cite[Lemma~2]{Alikhanov2015}.

\begin{lemma}
	\label{lem:6}
	Let $W^n\in L^2(\Omega)$ and $W^{n+\sigma}=\sigma W^{n+1}+(1-\sigma)W^n$. Then
	\begin{equation}
		2\big((\Ddisc_\tau^\beta W)^{n+\sigma},W^{n+\sigma}\big)_{L^2(\Omega)}
		\ge(\Ddisc_\tau^\beta\norm[L^2(\Omega)]{W}^2)^{n+\sigma}.
		\label{eq:36}
	\end{equation}
	Consequently,
	\begin{equation}
		\sum_{n=0}^{m-1}\tau
		\big((\Ddisc_\tau^\beta W)^{n+\sigma},W^{n+\sigma}\big)_{L^2(\Omega)}
		\ge-C_{\beta,T}\norm[L^2(\Omega)]{W^0}^2.
		\label{eq:37}
	\end{equation}
	For $W^0=0$, the left-hand side is nonnegative.
\end{lemma}
\begin{proof}
	Write $(\cdot,\cdot)=(\cdot,\cdot)_{L^2(\Omega)}$ and
	$\|\cdot\|=\|\cdot\|_{L^2(\Omega)}$.
	The scalar energy inequality in
	\cite[Lemma~1 and Corollary~1]{Alikhanov2015},
	with the weight properties in
	\cite[Lemma~4]{Alikhanov2015}, extends to $L^2(\Omega)$
	by polarization and proves \eqref{eq:36} for $n\ge1$.
	For $n=0$, it follows from
	\[
	2(W^1-W^0,W^\sigma)
	-\|W^1\|^2+\|W^0\|^2
	=(2\sigma-1)\|W^1-W^0\|^2\ge0.
	\]
	
	For the proof of \eqref{eq:37}, set $\nu=1-\beta=2\sigma-1$
	and define
	\[
	\mathcal G_nU
	=\sum_{k=0}^nC_{n-k}^{(n+1)}(U^{k+1}-U^k),
	\qquad
	\mathcal F_m(U)
	=\sum_{n=0}^{m-1}(\mathcal G_nU,U^{n+\sigma}).
	\]
	For $m=1$, completing the square gives
	\[
	\mathcal F_1(W)
	=A_0\left[\sigma\|W^1-W^0\|^2+(W^1-W^0,W^0)\right]
	\ge-\frac{A_0}{4\sigma}\|W^0\|^2,
	\]
	where $A_0=\sigma^\nu$.
	We may therefore assume $m\ge2$ and introduce
	\[
	\bar c_j=A_j+B_{j+1}-B_j\quad(j\ge0),
	\qquad
	e_0=A_0,\qquad e_j=A_j-B_j\quad(j\ge1),
	\]
	with $B_0=0$. The weight properties give
	\[
	\bar c_0>\bar c_1>\cdots>0,\qquad
	\bar c_{j-1}>e_j>0,\qquad
	(2\sigma-1)\bar c_0>\sigma\bar c_1.
	\]
	
	For $U^0=0$, let $\bm U=(U^1,\ldots,U^m)^\top$.
	Introduce matrices $D,P\in\mathbb R^{m\times m}$ such that
	\[
	D\bm U=(\mathcal G_0U,\ldots,\mathcal G_{m-1}U)^\top,
	\qquad
	P\bm U=(U^\sigma,\ldots,U^{m-1+\sigma})^\top.
	\]
	With $S=(D^\top P+P^\top D)/2$, we have
	\[
	\mathcal F_m(U)=\sum_{i,j=1}^mS_{ij}(U^i,U^j).
	\]
	The matrix $D$ is lower triangular with nonpositive
	off-diagonal entries. For $j<k$, multiplication gives
	\[
	2S_{jk}
	=\sigma D_{kj}
	+(1-\sigma)\mathbf1_{\{k<m\}}D_{k+1,j}
	+(1-\sigma)\mathbf1_{\{k=j+1\}}D_{kk},
	\]
	where $\mathbf1_{\{\cdot\}}$ is the indicator.
	Thus $S_{jk}\le0$ for $k\ge j+2$.
	For $k=j+1$, using
	$D_{kj}\le\bar c_1-\bar c_0$ and $D_{kk}=\bar c_0$ gives
	\[
	2S_{j,j+1}
	\le\sigma\bar c_1-(2\sigma-1)\bar c_0<0.
	\]
	
	The row sums $\eta_i=\sum_{j=1}^mS_{ij}$ satisfy
	\[
	\begin{aligned}
		2\eta_1
		&=2\sigma A_0+(1-\sigma)e_1-H_{m-1},\\
		2\eta_i
		&=\bar c_{m-i}+\sigma e_{i-1}+(1-\sigma)e_i,
		&&2\le i<m,\\
		2\eta_m&=\bar c_0+\sigma e_{m-1},
	\end{aligned}
	\]
	where
	\[
	H_N
	=\sum_{\ell=1}^N(\bar c_{\ell-1}-e_\ell)
	=A_0-A_N+\sum_{\ell=1}^{N-1}B_\ell+2B_N.
	\]
	The trapezoidal remainder for $s^\nu$ gives
	\[
	0<B_\ell\le\frac{\nu}{8}
	\left[
	(\ell-1+\sigma)^{\nu-1}
	-(\ell+\sigma)^{\nu-1}
	\right].
	\]
	Summation yields
	\[
	\sum_{\ell=1}^{\infty}B_\ell
	\le\frac{\nu}{8}\sigma^{\nu-1},
	\qquad
	H_N<A_0+\frac{\nu}{4}\sigma^{\nu-1}.
	\]
	Since $\sigma>1/2$ and $A_0=\sigma^\nu$, it follows that
	\[
	\eta_1
	\ge\frac{\nu A_0}{2}
	\left(1-\frac{1}{4\sigma}\right)
	\ge\frac{\nu A_0}{4}>0.
	\]
	The remaining row sums are positive by their expressions above.
	Consequently,
	\[
	\mathcal F_m(U)
	=\sum_{i=1}^m\eta_i\|U^i\|^2
	+\sum_{i<j}(-S_{ij})\|U^i-U^j\|^2
	\ge\sum_{i=1}^m\eta_i\|U^i\|^2.
	\]
	
	For arbitrary $W^0$, set $U^j=W^j-W^0$.
	Then $\mathcal G_nU=\mathcal G_nW$ and
	$W^{n+\sigma}=U^{n+\sigma}+W^0$.
	Telescoping gives
	\[
	\mathcal F_m(W)
	=\mathcal F_m(U)+\sum_{i=1}^mh_i(U^i,W^0),
	\qquad
	h_1=A_0-H_{m-1},\quad
	h_i=\bar c_{m-i}\quad(2\le i\le m).
	\]
	The preceding bounds imply
	\[
	|h_1|\le A_0,\qquad
	\eta_1\ge\frac{\nu A_0}{4},\qquad
	\eta_i\ge\frac{h_i}{2}\quad(2\le i\le m),
	\]
	and
	\[
	\sum_{i=2}^mh_i
	=(m-2+\sigma)^\nu+B_{m-1}
	\le C_\beta m^\nu.
	\]
	Completing the square therefore yields
	\[
	\mathcal F_m(W)
	\ge-\frac14\sum_{i=1}^m
	\frac{h_i^2}{\eta_i}\|W^0\|^2
	\ge-C_\beta m^{1-\beta}\|W^0\|^2.
	\]
	Multiplication by $\tau^{1-\beta}/\Gamma(2-\beta)$,
	together with $m\tau\le T$, proves \eqref{eq:37}.
\end{proof}
We recast the scheme in shifted-energy form by setting $\vartheta=\sigma-\frac{1}{2}$ and defining
\begin{equation*}
	\widetilde U_h^n=U_h^n+\vartheta\tau V_h^n,\qquad
	\widetilde V_h^n=V_h^n+\vartheta\tau A_h^n.
\end{equation*}
Equations \eqref{eq:31}--\eqref{eq:32} imply the exact identities
\begin{equation}
	\begin{aligned}
		U_h^{n+\sigma}
		&=\frac{\widetilde U_h^{n+1}+\widetilde U_h^n}{2},
		&\frac{\widetilde U_h^{n+1}-\widetilde U_h^n}{\tau}
		&=V_h^{n+\sigma},\\
		V_h^{n+\sigma}
		&=\frac{\widetilde V_h^{n+1}+\widetilde V_h^n}{2},
		&\frac{\widetilde V_h^{n+1}-\widetilde V_h^n}{\tau}
		&=A_h^{n+\sigma}.
	\end{aligned}
	\label{eq:38}
\end{equation}
Accordingly, define the shifted discrete energy by
\begin{equation*}
	\mathcal E_h^n
	=\frac{1}{2}m(\widetilde V_h^n,\widetilde V_h^n)
	+\frac{1}{2}a_h(\widetilde U_h^n,\widetilde U_h^n).
\end{equation*}

\begin{theorem}
	\label{thm:7}
	For every $1\le m\le N$,
	\begin{equation}
		\max_{0\le n\le m}\mathcal E_h^n
		\le C\left[
		\mathcal E_h^0+\muh C_{\beta,T}\norm[L^2(\Omega)]{V_h^0}^2
		+\sum_{n=0}^{m-1}\tau\norm[L^2(\Omega)]{f^{n+\sigma}}^2
		\right].
		\label{eq:39}
	\end{equation}

\end{theorem}

\begin{proof}
	By choosing $\chi_h = V_h^{n+\sigma}$ in \eqref{eq:33}, multiplying by $\tau$, and using \eqref{eq:38} and polarization, we obtain
	\begin{equation*}
		\mathcal E_h^{n+1}-\mathcal E_h^n
		+\tau d((\Ddisc_\tau^\beta V_h)^{n+\sigma},V_h^{n+\sigma})
		=\tau(f^{n+\sigma},V_h^{n+\sigma}).
	\end{equation*}
	For any $1\le r\le m$, summation and Lemma~\ref{lem:6}, with $d(w,z)=\muh(w,z)$, give
	\[
	\mathcal E_h^r
	\le K_0+\sum_{n=0}^{r-1}\tau
	\norm[L^2(\Omega)]{f^{n+\sigma}}\norm[L^2(\Omega)]{V_h^{n+\sigma}},
	\qquad
	K_0=\mathcal E_h^0+\muh C_{\beta,T}\norm[L^2(\Omega)]{V_h^0}^2.
	\]
	Moreover, \eqref{eq:38} and $\rhoh\ge\rho_*$ imply
	\[
	\norm[L^2(\Omega)]{V_h^{n+\sigma}}
	\le C\big((\mathcal E_h^{n+1})^{1/2}
	+(\mathcal E_h^n)^{1/2}\big).
	\]
	Choose the terminal index of the sum so that
	$M_m=\max_{0\le j\le m}(\mathcal E_h^j)^{1/2}$ is attained there. The Cauchy--Schwarz inequality
	in time then yields
	\[
	M_m^2\le K_0+C\sqrt{T}\,M_m
	\left(\sum_{n=0}^{m-1}\tau\norm[L^2(\Omega)]{f^{n+\sigma}}^2\right)^{1/2}.
	\]
    Young's inequality proves \eqref{eq:39}.
\end{proof}

\subsection{Error analysis for microscopic cell functions and macroscopic homogenized material parameters}

Let $Y_1$ and $Y_2$ denote the matrix and inclusion subdomains
of the PUC $Y$, respectively. Let $p_Y\ge 1$ be the polynomial degree
of the cell finite element space. For $m=p_Y+1$, define the
piecewise Sobolev space by
\[
H_{\mathrm{pw}}^m(Y)
=\left\{v\in L^2(Y):
v|_{Y_r}\in H^m(Y_r),\ r=1,2\right\}.
\]

\begin{lemma}
	\label{lem:8}
	Suppose that $N^k,F,G,N^{k\ell}\in H_{\mathrm{pw}}^{p_Y+1}(Y)$ for $1\le k,\ell\le d$. Then
	\begin{align}
		&\max_k\norm[L^2(Y)]{\nabla_{\by}(N^k-N_{h_Y}^k)}
		+\norm[L^2(Y)]{\nabla_{\by}(F-F_{h_Y})}
		+\norm[L^2(Y)]{\nabla_{\by}(G-G_{h_Y})}\notag\\
		&\qquad+\max_{k,\ell}\norm[L^2(Y)]{\nabla_{\by}(N^{k\ell}-N_{h_Y}^{k\ell})}
		\le Ch_Y^{p_Y}.
	\end{align}

\end{lemma}
\begin{proof}
	These estimates follow directly from standard finite element approximation theory.
\end{proof}
Let $\mathcal{L}(\mathbb{R}^d)$ denote the space of linear operators on $\mathbb{R}^d$. For a matrix $B\in\mathbb{R}^{d\times d}$, define its operator norm by
\[
\|B\|_{\mathcal{L}(\mathbb{R}^d)}
:=\sup_{\boldsymbol{\xi}\in\mathbb{R}^d\setminus\{0\}}
\frac{|B\boldsymbol{\xi}|}{|\boldsymbol{\xi}|},
\]
where $|\cdot|$ denotes the Euclidean norm.
\begin{lemma}
	\label{lem:9}
	Under the assumptions of Lemma~\ref{lem:8}, the macroscopic homogenized  material  parameter in \eqref{eq:28} satisfies
	\begin{equation}
		\norm[\mathcal L(\R^d)]{\Ah_{h_Y}-\Ah}\le Ch_Y^{2p_Y}.
	\end{equation}
\end{lemma}
\begin{proof}
    Set $e^k=N^k-N_{h_Y}^k$. As in \cite{DongCao2009},
    the weak cell equation, symmetry of $A$, and Galerkin
    orthogonality yield
    \[
    \begin{aligned}
    \left|\widehat a_{k\ell,h_Y}-\widehat a_{k\ell}\right|
    &=\frac{1}{|Y|}\left|\mathfrak a_Y(e^k,e^\ell)\right|\\
    &\le C\norm[L^2(Y)]{\nabla_{\by}e^k}
    \norm[L^2(Y)]{\nabla_{\by}e^\ell}
    \le Ch_Y^{2p_Y},
    \end{aligned}
    \]
    where the last inequality follows from Lemma~\ref{lem:8}.
    Combining the componentwise bounds proves the result.
\end{proof}

\subsection{Error analysis for fully discrete method at the macro-scale}

The macroscopic homogenized material parameter \(\Ah_{h_Y}\) used in the computation is obtained from finite element approximations of the first-order auxiliary cell functions. Consequently, the fully discrete scheme discretizes a modified macroscopic homogenized equation with $\Ah$ replaced by $\Ah_{h_Y}$. Let $\bar u_0$ denote the exact solution of the resulting modified homogenized problem
\begin{equation}
	\left\{\begin{aligned}
		\rhoh\bar u_{0,tt}+\muh\CD{\alpha}\bar u_0
		-\nabla\cdot(\Ah_{h_Y}\nabla\bar u_0)&=f,
		&&\text{in }\Om\times(0,T],\\
		\bar u_0&=\widehat u(\bx),&&\text{on }\partial\Om\times[0,T],\\
		\bar u_0(0)=g,\qquad\bar u_{0,t}(0)&=\psi,&&\text{in }\Om.
	\end{aligned}\right.
	\label{eq:45}
\end{equation}
\begin{lemma}
	\label{lem:10}
	Let $\bar u_0$ denote the exact solution of the modified homogenized equation \eqref{eq:45}. Under the assumptions of Lemma~\ref{lem:9}, suppose that $u_0\in C^1([0,T];H^1(\Om))$. Then
	\begin{equation}
		\sup_{0\le t\le T}\left(
		\norm[L^2(\Omega)]{\bar u_{0,t}(t)-u_{0,t}(t)}
		+\norm[H^1(\Omega)]{\bar u_0(t)-u_0(t)}\right)
		\le C_T h_Y^{2p_Y},
		\label{eq:46}
	\end{equation}
	where $C_T$ is independent of $h_Y$.
\end{lemma}
\begin{proof}
	Denote $w=\bar u_0-u_0$ and $\Delta A=\Ah_{h_Y}-\Ah$. Subtracting the macroscopic homogenized equation \eqref{eq:13} from the modified homogenized equation \eqref{eq:45} gives
	\begin{equation}
		\rhoh w_{tt}+\muh\CD{\beta}w_t-\nabla\cdot(\Ah_{h_Y}\nabla w) =\nabla\cdot(\Delta A\nabla u_0),
		\label{eq:47}
	\end{equation}
	with $w=0$ on $\partial\Om\times[0,T]$ and homogeneous initial conditions $w(\bx,0)=w_t(\bx,0)=0$. 
	Testing \eqref{eq:47} with $w_t$ and applying Green's formula give
	\begin{align}
		&\frac{1}{2}\frac{\dd}{\dd t}\left[\int_{\Om}\rhoh |w_t|^2\dd\bx\right]
		+\frac{1}{2}\frac{\dd}{\dd t}\left[\int_{\Om}\Ah_{h_Y}\nabla w\cdot\nabla w\dd\bx\right]
		+\int_{\Om}\muh\CD{\beta}w_t w_t\dd\bx\notag\\
		&\qquad=-\int_{\Om}\Delta A\nabla u_0\cdot\nabla w_t\dd\bx.
		\label{eq:48}
	\end{align}
	Furthermore, integrating the right-hand side of \eqref{eq:48} by parts in time yields
	\begin{align}
		-\int_0^t\int_{\Om}\Delta A\nabla u_0\cdot\nabla w_t\dd\bx\dd s
		=&-\int_{\Om}\Delta A\nabla u_0(t)\cdot\nabla w(t)\dd\bx\notag\\
		&+\int_0^t\int_{\Om}\Delta A\nabla u_{0,t}(s)\cdot\nabla w(s)\dd\bx\dd s.
		\label{eq:49}
	\end{align}
	Integrate \eqref{eq:48} from $0$ to $t$ and use \eqref{eq:49}. Lemmas~\ref{lem:9} and~\ref{lem:1}, the Cauchy--Schwarz inequality, and Young's inequality then give
	\begin{align}
		&\frac{1}{2}\int_{\Om}\rhoh |w_t(t)|^2\dd\bx
		+\frac{1}{2}\int_{\Om}\Ah_{h_Y}\nabla w(t)\cdot\nabla w(t)\dd\bx
		+\frac{1}{2}\int_0^t\CD{\beta}\norm[L^2(\Omega)]{\sqrt{\muh}w_t(s)}^2\dd s\notag\\
		&\quad\le
		\norm[\mathcal L(\R^d)]{\Delta A}\norm[L^2(\Omega)]{\nabla u_0(t)}\norm[L^2(\Omega)]{\nabla w(t)}\notag\\
		&\qquad+\int_0^t\norm[\mathcal L(\R^d)]{\Delta A}\norm[L^2(\Omega)]{\nabla u_{0,t}(s)}\norm[L^2(\Omega)]{\nabla w(s)}\dd s\notag\\
		&\quad\le
		\frac{1}{4}\int_{\Om}\Ah_{h_Y}\nabla w(t)\cdot\nabla w(t)\dd\bx
		+C\norm[\mathcal L(\R^d)]{\Delta A}^2\norm[H^1(\Omega)]{u_0(t)}^2\notag\\
		&\qquad+\int_0^t C\norm[\mathcal L(\R^d)]{\Delta A}^2\norm[H^1(\Omega)]{u_{0,t}(s)}^2\dd s
		+\int_0^t C\norm[H^1(\Omega)]{w(s)}^2\dd s\notag\\
		&\quad\le
		\frac{1}{4}\int_{\Om}\Ah_{h_Y}\nabla w(t)\cdot\nabla w(t)\dd\bx
		+C_T h_Y^{4p_Y}\left(\norm[L^\infty(0,T;H^1(\Omega))]{u_0}^2+\norm[L^2(0,T;H^1(\Omega))]{u_{0,t}}^2\right)\notag\\
		&\qquad+\int_0^t C\norm[H^1(\Omega)]{w(s)}^2\dd s.
		\label{eq:50}
	\end{align}
	Set $\Upsilon(t)=\norm[L^2(\Omega)]{w_t(t)}^2+\norm[H^1(\Omega)]{w(t)}^2$. Then \eqref{eq:50} gives $\Upsilon(t)\le C_T h_Y^{4p_Y}+C\int_0^t\Upsilon(s)\dd s$. Consequently, Gronwall's inequality yields, for every $t\in[0,T]$,
	\begin{equation}
		\norm[L^2(\Omega)]{w_t(t)}^2+\norm[H^1(\Omega)]{w(t)}^2 \le C_T h_Y^{4p_Y}.
	\end{equation}
	Taking square roots proves \eqref{eq:46}.
\end{proof}
\begin{lemma}
	\label{lem:11}
	Let $\bar u_0$ be the solution of
	\eqref{eq:45}. Assume that
	$f\in C^2([0,T];L^2(\Om))$ and
	\begin{equation}
		\bar u_0
		\in C^4([0,T];H^1(\Om))
		\cap C^2([0,T];H^{p+1}(\Om)),
		\label{eq:52}
	\end{equation}
	with the corresponding norms bounded independently of $h_Y$.
	Suppose also that
	\begin{equation}
		\bar u_{0,tt}(0)=0,
		\qquad
		\bar u_{0,ttt}(0)=0.
		\label{eq:53}
	\end{equation}
	Under the stated spatial discretization assumptions,
	the scheme
	\eqref{eq:31}--\eqref{eq:33},
	initialized by $U_h^0=R_h^Dg$, $V_h^0=P_h\psi$, and
	\eqref{eq:34}, satisfies
	\begin{equation}
		\max_{0\le n\le N}
		\left(
		\norm[L^2(\Omega)]{V_h^n-\bar u_{0,t}(t_n)}
		+
		\norm[H^1(\Omega)]{U_h^n-\bar u_0(t_n)}
		\right)
		\le C_T(h^p+\tau^2).
		\label{eq:54}
	\end{equation}
	The same estimate holds at $t_{n+\sigma}$ for the
	linearly interpolated displacement and velocity.
\end{lemma}
\begin{proof}
Write $u^n=\bar u_0(t_n)$, $v^n=\bar u_{0,t}(t_n)$, and
$a^n=\bar u_{0,tt}(t_n)$.
Define comparison sequences satisfying the Newmark relations by
\begin{align*}
\overline A^n&=R_ha^n,\qquad
\overline V^0=R_h\psi,\qquad \overline U^0=R_h^Dg,\\
\overline V^{n+1}-\overline V^n
&=\frac{\tau}{2}(\overline A^{n+1}+\overline A^n),\\
\overline U^{n+1}-\overline U^n
&=\frac{\tau}{2}(\overline V^{n+1}+\overline V^n).
\end{align*}
Set $e_Z^n=Z_h^n-\overline Z^n$ for $Z=U,V,A$, and
$q_U^n=\overline U^n-R_h^Du^n$, $q_V^n=\overline V^n-R_hv^n$.
Then
\begin{equation}
\begin{aligned}
U_h^n-u^n&=e_U^n+q_U^n+(R_h^Du^n-u^n),\\
V_h^n-v^n&=e_V^n+q_V^n+(R_hv^n-v^n).
\end{aligned}
\label{eq:55}
\end{equation}
For $w\in H_0^1(\Om)\cap H^{p+1}(\Om)$, the standard Ritz estimates
\cite{BrennerScott2008} give
\begin{equation}
\norm[H^1(\Om)]{w-R_hw}
+h^{-1}\norm[L^2(\Om)]{w-R_hw}
\le Ch^p\norm[H^{p+1}(\Om)]{w}.
\label{eq:62}
\end{equation}
Since $\operatorname{Tr}(I_hu^n)=\widehat u_h$,
Ritz orthogonality and interpolation also yield
$\norm[H^1(\Om)]{R_h^Du^n-u^n}
\le C\norm[H^1(\Om)]{I_hu^n-u^n}\le Ch^p$.

The trapezoidal remainders in the comparison recurrences, together with
$q_U^0=q_V^0=0$ and \eqref{eq:52}, give
\begin{equation}
\begin{aligned}
&\max_{0\le n\le N}
\left(\norm[H^1(\Om)]{q_U^n}+\norm[H^1(\Om)]{q_V^n}\right)
+\sum_{n=0}^{N-1}\norm[H^1(\Om)]{q_U^{n+1}-q_U^n}
\le C_T\tau^2,\\
&\max_{0\le n<N}\norm[H^1(\Om)]{q_V^{n+1}-q_V^n}
\le C\tau^3.
\end{aligned}
\label{eq:61}
\end{equation}
Positivity of the convolution weights and
$\sum_{j=0}^n C_j^{(n+1)}=(n+\sigma)^{1-\beta}$ therefore imply
\[
\norm[L^2(\Om)]{(\Ddisc_\tau^\beta q_V)^{n+\sigma}}
\le C\tau^{3-\beta}(n+\sigma)^{1-\beta}\le C_T\tau^2.
\]
Let $g_\beta=\CD{\beta}\bar u_{0,t}$ and
$w_I^{n+\sigma}=\sigma w(t_{n+1})+(1-\sigma)w(t_n)$.
Let $I$ denote the identity operator. Define the residual by
\begin{align*}
r^{n+\sigma}
={}&\rhoh(R_h-I)a_I^{n+\sigma}
+\muh\left[R_h(\Ddisc_\tau^\beta v)^{n+\sigma}
-(g_\beta)_I^{n+\sigma}\right]\\
&+\muh(\Ddisc_\tau^\beta q_V)^{n+\sigma}
+f_I^{n+\sigma}-f(t_{n+\sigma}).
\end{align*}
The regularity and compatibility conditions imply
$g_\beta\in C^2([0,T];H_0^1(\Om))\cap C([0,T];H^{p+1}(\Om))$.
Thus \eqref{eq:35}, \eqref{eq:62}, and the interpolation remainder give
\[
\max_{0\le n<N}\norm[L^2(\Om)]{r^{n+\sigma}}
\le C_T(h^{p+1}+\tau^2).
\]
Subtracting the projected, endpoint-interpolated equation \eqref{eq:45}
from \eqref{eq:33} yields, for every $\chi_h\in V_h$,
\begin{equation}
\begin{aligned}
&m(e_A^{n+\sigma},\chi_h)
+d((\Ddisc_\tau^\beta e_V)^{n+\sigma},\chi_h)
+a_h(e_U^{n+\sigma},\chi_h)\\
&\qquad=-(r^{n+\sigma},\chi_h)-a_h(q_U^{n+\sigma},\chi_h).
\end{aligned}
\label{eq:57}
\end{equation}
With $\vartheta=\sigma-\tfrac12$, set
\[
\begin{gathered}
\widetilde e_U^n=e_U^n+\vartheta\tau e_V^n,\qquad
\widetilde e_V^n=e_V^n+\vartheta\tau e_A^n,\\
\mathcal E_e^n=\tfrac12 m(\widetilde e_V^n,\widetilde e_V^n)
+\tfrac12 a_h(\widetilde e_U^n,\widetilde e_U^n).
\end{gathered}
\]
Both sequences satisfy the Newmark relations. We therefore test \eqref{eq:57}
with $\tau e_V^{n+\sigma}$ and use \eqref{eq:38} to obtain
\begin{equation}
\begin{aligned}
\mathcal E_e^{n+1}-\mathcal E_e^n
&+\tau d((\Ddisc_\tau^\beta e_V)^{n+\sigma},e_V^{n+\sigma})\\
&=-\tau(r^{n+\sigma},e_V^{n+\sigma})
-a_h(q_U^{n+\sigma},\widetilde e_U^{n+1}-\widetilde e_U^n).
\end{aligned}
\label{eq:58}
\end{equation}
The initialization and \eqref{eq:53} give
$e_U^0=e_A^0=0$ and $e_V^0=P_h\psi-R_h\psi$.
Projection estimates and the $H^1$ stability of $P_h$
\cite{BrennerScott2008,CrouzeixThomee1987} imply
\[
(\mathcal E_e^0)^{1/2}+\norm[L^2(\Om)]{e_V^0}
\le C(h^{p+1}+\tau h^p)\le C(h^{p+1}+\tau^2).
\]
Let $M_N=\max_{0\le n\le N}(\mathcal E_e^n)^{1/2}$.
Discrete summation by parts and \eqref{eq:61} yield, for $1\le r\le N$,
\[
\left|\sum_{n=0}^{r-1}
a_h(q_U^{n+\sigma},\widetilde e_U^{n+1}-\widetilde e_U^n)\right|
\le C_T\tau^2 M_N.
\]
Summing \eqref{eq:58} and applying Lemma~\ref{lem:6}, as in
Theorem~\ref{thm:7}, gives
\[
\begin{gathered}
M_N^2\le C_T(h^{p+1}+\tau^2)^2
+C_T(h^{p+1}+\tau^2)M_N,\\
M_N\le C_T(h^{p+1}+\tau^2).
\end{gathered}
\]
For $Z=U,V$, the identity
$\sigma e_Z^{n+1}+(1-\sigma)e_Z^n
=(\widetilde e_Z^{n+1}+\widetilde e_Z^n)/2$
and $(1-\sigma)/\sigma<1$ imply
\[
\max_n\norm[X_Z]{e_Z^n}
\le\norm[X_Z]{e_Z^0}
+\frac{1}{2\sigma-1}\max_n\norm[X_Z]{\widetilde e_Z^n},
\qquad X_U=H^1(\Om),\quad X_V=L^2(\Om).
\]
Hence
\begin{equation}
\max_{0\le n\le N}
\left(\norm[H^1(\Om)]{e_U^n}+\norm[L^2(\Om)]{e_V^n}\right)
\le C_T(h^{p+1}+\tau^2).
\label{eq:60}
\end{equation}
Combining \eqref{eq:55}, \eqref{eq:61}, \eqref{eq:60}, and the Ritz
estimates proves \eqref{eq:54}. Linear interpolation and its
$O(\tau^2)$ remainder give the shifted-time estimate.
\end{proof}

Combining the coefficient-perturbation estimate with the fully discrete error estimate yields the following total error bound for the macroscopic homogenized solution.

\begin{theorem}
	Let \(U_h^n\) and \(V_h^n\) be the fully discrete approximations to the displacement and velocity obtained by the scheme \eqref{eq:31}--\eqref{eq:33} for the modified homogenized problem \eqref{eq:45}, and let \(u_0\) be the exact solution of the homogenized problem \eqref{eq:13}. Under the hypotheses of Lemmas~\ref{lem:8}--\ref{lem:11}, the following estimate holds:
	\[
	\begin{aligned}
		&\max_{0\le n\le N}
		\left\|U_h^n-u_0(t_n)\right\|_{H^1(\Omega)}+
		\max_{0\le n\le N}
		\left\|V_h^n-u_{0,t}(t_n)\right\|_{L^2(\Omega)}\\
		\le& C_T\left(h_Y^{2p_Y}+h^p+\tau^2\right),
	\end{aligned}
	\]
	where \(C_T\) is independent of \(h_Y\), \(h\), and \(\tau\), under the uniform bounds assumed in these lemmas. In particular, for \(p_Y=1\) and \(p=2\), the right-hand side becomes
	\[
	C_T\left(h_Y^2+h^2+\tau^2\right).
	\]
\end{theorem}
\begin{proof}
    The result follows from Lemmas~\ref{lem:10}
    and~\ref{lem:11}, and the triangle inequality.
\end{proof}

\begingroup
\renewcommand{\topfraction}{0.95}
\renewcommand{\bottomfraction}{0.95}
\renewcommand{\textfraction}{0.05}
\renewcommand{\floatpagefraction}{0.75}
\setcounter{totalnumber}{5}
\setlength{\textfloatsep}{10pt plus 2pt minus 2pt}
\setlength{\floatsep}{10pt plus 2pt minus 2pt}
\setlength{\intextsep}{10pt plus 2pt minus 2pt}
\makeatletter
\setlength{\@fptop}{0pt}
\setlength{\@fpsep}{10pt}
\setlength{\@fpbot}{0pt plus 1fil}
\makeatother
\section{Numerical examples and results}
\label{sec:6}
All computations were performed in FreeFEM++ on the same HPS desktop workstation, equipped with an Intel(R) Xeon(R) Gold 6146 CPU at $3.20$ GHz and $96.0$ GB of memory. All examples use $\alpha=3/2$ and $\sigma=3/4$, with zero initial displacement and velocity:
\[
u^\eps(\bx,0)=u_t^\eps(\bx,0)=0,\qquad \bx\in\Om.
\]
For $0<t\le T$, Examples~1--3 impose $u^\eps=0$ on $\partial\Om$. In Examples~4 and~5, fixed supports are imposed on the four
lateral faces, while the top and bottom faces are traction-free. Tables~\ref{tab:materials} and~\ref{tab:cost} collect the material parameters, mesh statistics, and CPU times.

The preceding error estimates concern Dirichlet boundary conditions; Examples~4 and~5 provide numerical tests with mixed boundary conditions. In Examples~1--4, DNS solves the original multi-scale problem on a fine mesh with the same loading, initial and boundary data, and time step as the HOMS computation. At $t_{n+\sigma}$, let $u_h^{[j],n+\sigma}$ denote the homogenized, LOMS, and HOMS displacement fields for $j=0,1,2$, respectively. The DNS reference solution is interpolated from endpoint states as
\[
u_{\rm DNS}^{\eps,n+\sigma}
=\sigma u_{\rm DNS}^{\eps,n+1}+(1-\sigma)u_{\rm DNS}^{\eps,n}.
\]
The relative $L^2$-errors and $H^1$-seminorm errors are
\begin{align*}
E_{L^2}^{[j]}(t_{n+\sigma})
&=\frac{\norm[L^2(\Omega)]{u_h^{[j],n+\sigma}-u_{\rm DNS}^{\eps,n+\sigma}}}{\norm[L^2(\Omega)]{u_{\rm DNS}^{\eps,n+\sigma}}},\\
E_{H^1}^{[j]}(t_{n+\sigma})
&=\frac{\norm[L^2(\Omega)]{\nabla(u_h^{[j],n+\sigma}-u_{\rm DNS}^{\eps,n+\sigma})}}{\norm[L^2(\Omega)]{\nabla u_{\rm DNS}^{\eps,n+\sigma}}}.
\end{align*}
For these comparisons, the quantities $E_{L^2}^{[0]}(t)$, $E_{L^2}^{[1]}(t)$, and $E_{L^2}^{[2]}(t)$ denote the relative $L^2$-errors of the homogenized, LOMS, and HOMS solutions, respectively. Similarly, $E_{H^1}^{[0]}(t)$, $E_{H^1}^{[1]}(t)$, and $E_{H^1}^{[2]}(t)$ denote the corresponding relative $H^1$-seminorm errors. Each approximate solution is compared with the DNS reference at the same shifted time $t_{n+\sigma}$.

\begin{table}[!htbp]
\centering
\caption{Material parameters for all numerical examples. Entries are matrix/inclusion values. For $d=2$, density denotes areal mass density.}
\label{tab:materials}
\small
\setlength{\tabcolsep}{0pt}
\begin{tabular}{@{}>{\centering\arraybackslash}p{0.13\textwidth}>{\centering\arraybackslash}p{0.05\textwidth}>{\centering\arraybackslash}p{0.24\textwidth}>{\centering\arraybackslash}p{0.34\textwidth}>{\centering\arraybackslash}p{0.24\textwidth}@{}}
\toprule
Examples & $d$ & Stiffness coefficient & Fractional damping coefficient & Mass density\\
\midrule
1--2 & 2 & $6000/240$ & $100/1$ & $10/1$\\
3 & 2 & $1000/10$ & $50/0.5$ & $20/1$\\
4--5 & 3 & $1000/50$ & $50/0.5$ & $20/1$\\
\midrule
Units & 2 & $\mathrm{kg}/\mathrm{s}^{2}$ & $\mathrm{kg}/(\mathrm{cm}^{2}\cdot\mathrm{s}^{1/2})$ & $\mathrm{kg}/\mathrm{cm}^{2}$\\
 & 3 & $\mathrm{kg}/(\mathrm{cm}\cdot\mathrm{s}^{2})$ & $\mathrm{kg}/(\mathrm{cm}^{3}\cdot\mathrm{s}^{1/2})$ & $\mathrm{kg}/\mathrm{cm}^{3}$\\
\bottomrule
\end{tabular}
\end{table}

\begin{table}[!htbp]
\centering
\caption{Mesh statistics and CPU times.}
\label{tab:cost}
\small
\setlength{\tabcolsep}{5pt}
\begin{tabular}{@{}clrrr@{}}
\toprule
Example & Quantity & Multi-scale eqs. & Cell eqs. & Homogenized eqs.\\
\midrule
1 & FEM elements & $297\,344$ & $484$ & $3\,200$\\
  & FEM nodes & $149\,153$ & $263$ & $1\,681$\\
\addlinespace[2pt]
2 & FEM elements & $182\,330$ & $2\,704$ & $3\,200$\\
  & FEM nodes & $91\,250$ & $1\,419$ & $1\,681$\\
\addlinespace[2pt]
3 & FEM elements & $819\,200$ & $1\,800$ & $3\,200$\\
  & FEM nodes & $410\,881$ & $961$ & $1\,681$\\
\addlinespace[2pt]
4 & FEM elements & $663\,552$ & $3\,072$ & $6\,000$\\
  & FEM nodes & $117\,649$ & $729$ & $1\,331$\\
\addlinespace[2pt]
5 & FEM elements & $\approx3.93\times10^7$ & $3\,072$ & $76\,800$\\
  & FEM nodes & $\approx6.70\times10^6$ & $729$ & $15\,129$\\
\midrule
\multicolumn{2}{c}{Example} & DNS time (s) & \multicolumn{2}{c}{HOMS time (s)}\\
\midrule
\multicolumn{2}{c}{1} & $21521.7$ & \multicolumn{2}{c}{$2264.06$}\\
\multicolumn{2}{c}{2} & $11253.6$ & \multicolumn{2}{c}{$2117.05$}\\
\multicolumn{2}{c}{3} & $73676.6$ & \multicolumn{2}{c}{$29024.4$}\\
\multicolumn{2}{c}{4} & $7047.94$ & \multicolumn{2}{c}{$1428.09$}\\
\bottomrule
\end{tabular}
\end{table}

\subsection{Validation of the HOMS computational method}
\subsubsection{Example 1: The nonlocal-in-time fractional wave problem in 2D heterogeneous structures with circular inclusions}
\label{subsec:6.1}

The first medium contains circular inclusions in the square domain $\Omega=[0,1\,\mathrm{cm}]^2$, with $\eps=1/8$. Each inclusion has radius $0.03125\,\mathrm{cm}$. Figure~\ref{fig:1} shows the domain and the PUC $Y$. Table~\ref{tab:materials} lists the material parameters for Example~1. The wave response is simulated over $t\in[0,2]\,\mathrm{s}$ using the time step $\tau=0.01\,\mathrm{s}$.

\begin{figure}[!htbp]
	\centering
	\includegraphics[width=0.82\textwidth]{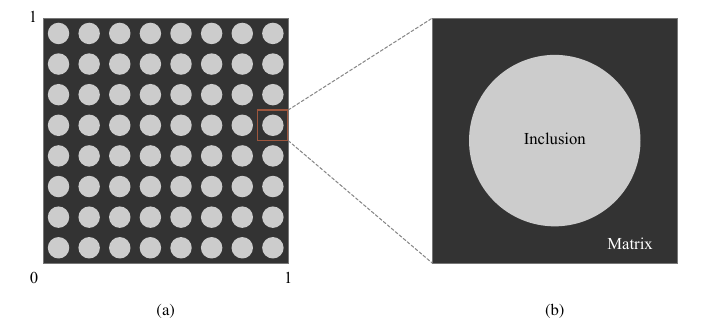}
	\caption{Example~1. (a) The 2D composite material $\Omega$; (b) reference PUC $Y$.}
	\label{fig:1}
\end{figure}

The source term is
\[
\begin{aligned}
f(\bx,t)
=100\sin\left(\pi x_1\right)
\sin\left(\pi x_2\right)
\left[1.5-\cos(4\pi t)\right]
\,\mathrm{kg}/(\mathrm{cm}\cdot\mathrm{s}^{2}),\quad 
\bx\in\Om,\ 0<t\le T.
\end{aligned}
\]
 Table~\ref{tab:cost} summarizes the mesh statistics and CPU times for Example~1. Figure~\ref{fig:2} compares the four wave solutions, and Figure~\ref{fig:3} shows the corresponding absolute errors, at $t_{40+\sigma}=0.4075\,\mathrm{s}$. Figure~\ref{fig:4} shows the time histories of the relative $L^2$-errors and relative $H^1$-seminorm errors.

\begin{figure}[!htbp]
	\centering
	\includegraphics[width=0.92\textwidth]{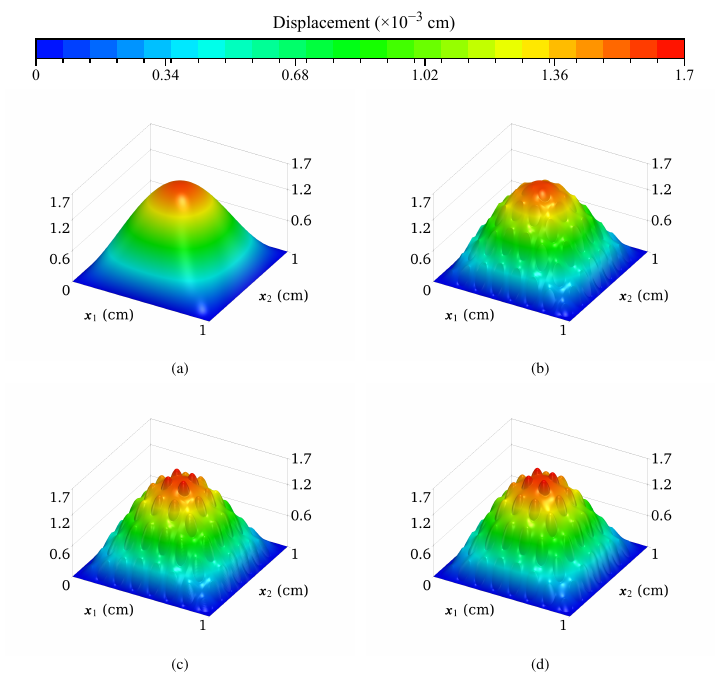}
	\caption{Example~1. Displacement fields: (a) $u_0$; (b) $u^{(1\varepsilon)}$; (c) $u^{(2\varepsilon)}$; (d) $u_{\rm DNS}^{\varepsilon}$.}
	\label{fig:2}

    \vspace{6pt}
\includegraphics[width=0.92\textwidth]{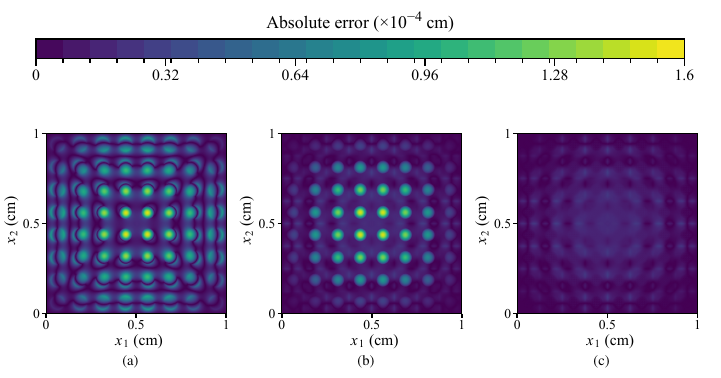}
\caption{Example~1. Absolute errors of displacement fields:
	(a) $|u_0-u_{\rm DNS}^{\varepsilon}|$;
	(b) $|u^{(1\varepsilon)}-u_{\rm DNS}^{\varepsilon}|$;
	(c) $|u^{(2\varepsilon)}-u_{\rm DNS}^{\varepsilon}|$.}
\label{fig:3}
\end{figure}

\begin{figure}[!htbp]
	\centering
	\includegraphics[width=0.96\textwidth]{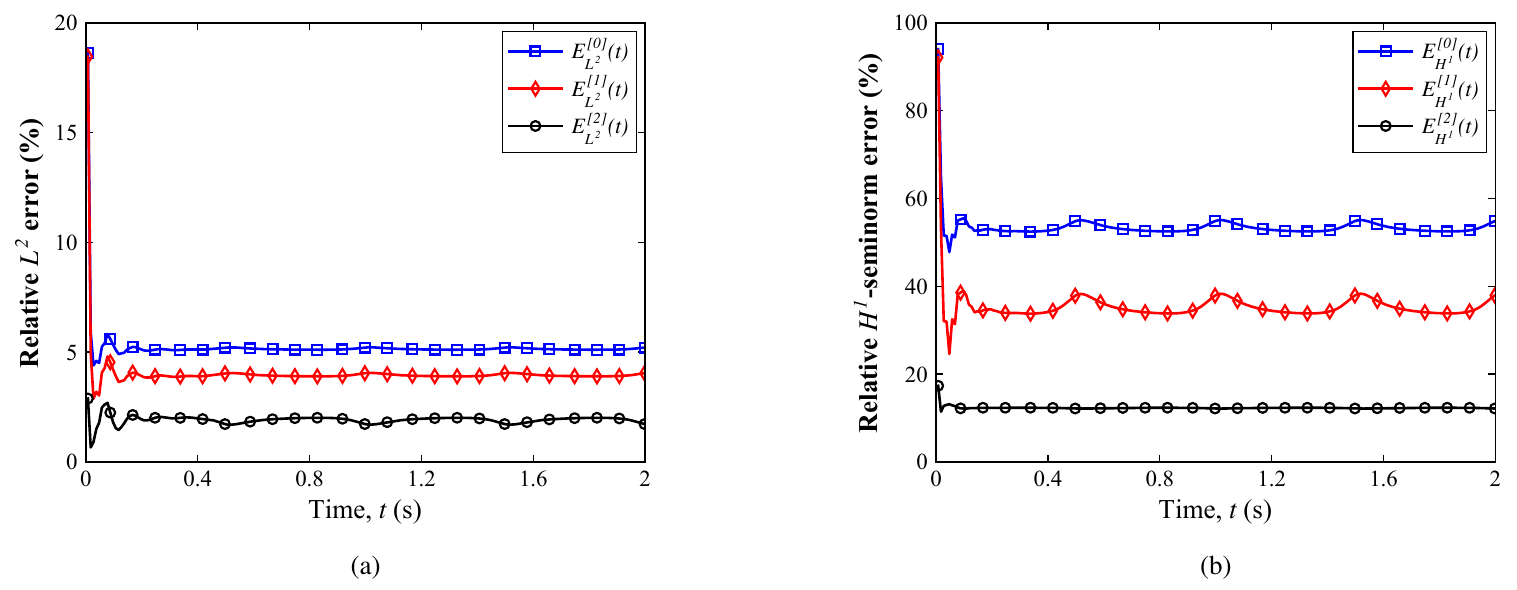}
	\caption{Example~1. Time evolution of relative errors: (a) $E_{L^2}^{[j]}(t)$; (b) $E_{H^1}^{[j]}(t)$, $j=0,1,2$.}
	\label{fig:4}
\end{figure}

Figure~\ref{fig:4} shows that the HOMS solution maintains the lowest relative errors over most of the displayed time interval. No numerical blow-up is observed up to $t=2.0\,\mathrm{s}$, supporting numerical robustness over the simulated interval. The HOMS computation achieves a single-run speedup of approximately $9.51$ over DNS. These results demonstrate the accuracy of the HOMS solution at substantially lower computational cost than DNS.

\subsubsection{Example 2: The nonlocal-in-time fractional wave problem in 2D heterogeneous structures with diamond-shaped inclusions }
\label{subsec:6.2}

With $\eps=1/8$, the second medium uses the same domain $\Omega$ and material parameters as Example~1. Each PUC contains a centered diamond-shaped inclusion with both diagonals of length $0.075\,\mathrm{cm}$. Figure~\ref{fig:5} shows the composite and its PUC.
The wave response is simulated over $t\in[0,1]\,\mathrm{s}$ with $\tau=0.005\,\mathrm{s}$.
\begin{figure}[!htbp]
	\centering
	\includegraphics[width=0.82\textwidth]{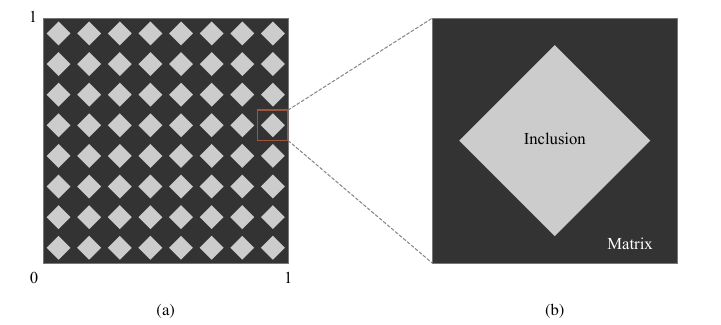}
	\caption{Example~2. (a) The 2D composite material $\Omega$; (b) reference PUC $Y$.}
	\label{fig:5}
\end{figure}
The source term is
\[
f(\bx,t)=600\,\mathrm{kg}/(\mathrm{cm}\cdot\mathrm{s}^{2}),
\qquad \bx\in\Omega,\quad 0<t\le T.
\]
Table~\ref{tab:cost} summarizes the mesh statistics and CPU times for Example~2. Figure~\ref{fig:6} compares the homogenized, LOMS, HOMS, and DNS solutions at $t_{25+\sigma}=0.12875\,\mathrm{s}$, while Figure~\ref{fig:7} shows the corresponding absolute errors. Figure~\ref{fig:8} presents the time evolution of the relative $L^2$-errors and relative $H^1$-seminorm errors.

The absolute errors are concentrated near the diamond-shaped interfaces and decrease with the order of the multi-scale asymptotic solution. The HOMS solution gives the lowest displacement and gradient errors over most of the displayed time interval. The $H^1$-seminorm comparison is particularly relevant to recovery of the local gradients at material interfaces. The DNS mesh contains approximately $56.98$ times as many elements as the homogenized mesh, so the macroscopic system is substantially smaller. The HOMS computation achieves a single-run speedup of approximately $5.32$ over DNS. These results support the accuracy and efficiency of the HOMS asymptotic solution for this diamond-inclusion composite under constant loading.

\begin{figure}[!htbp]
	\centering
	\includegraphics[width=0.92\textwidth]{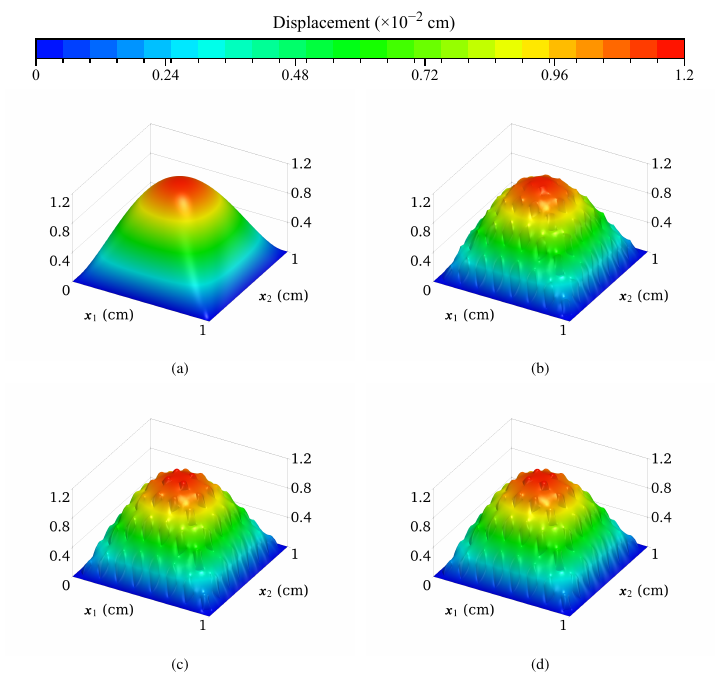}
	\caption{Example~2. Displacement fields: (a) $u_0$; (b) $u^{(1\varepsilon)}$; (c) $u^{(2\varepsilon)}$; (d) $u_{\rm DNS}^{\varepsilon}$.}
	\label{fig:6}

    \vspace{6pt}
\includegraphics[width=0.92\textwidth]{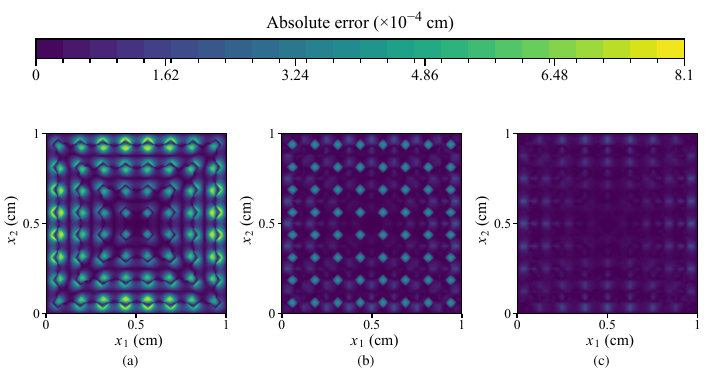}
\caption{Example~2. Absolute errors of displacement fields:
	(a) $|u_0-u_{\rm DNS}^{\varepsilon}|$;
	(b) $|u^{(1\varepsilon)}-u_{\rm DNS}^{\varepsilon}|$;
	(c) $|u^{(2\varepsilon)}-u_{\rm DNS}^{\varepsilon}|$.}
	\label{fig:7}
\end{figure}

\begin{figure}[!htbp]
	\centering
	\includegraphics[width=0.96\textwidth]{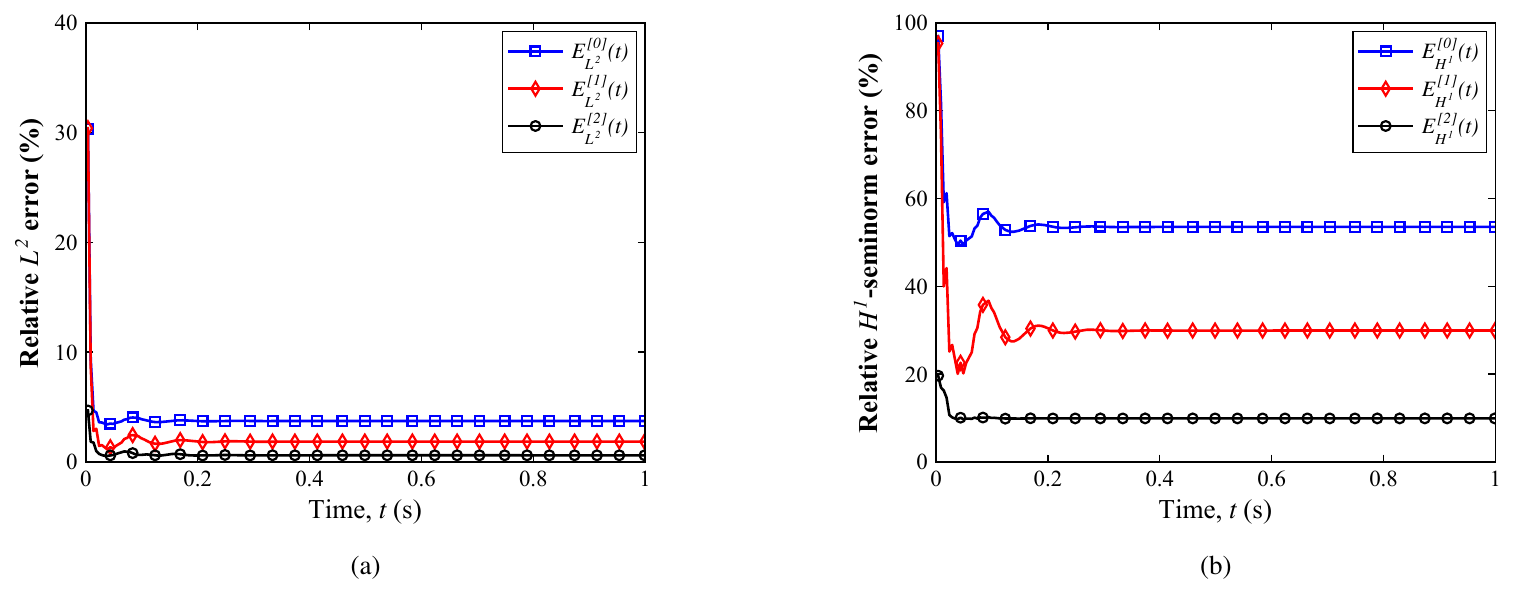}
	\caption{Example~2. Time evolution of relative errors: (a) $E_{L^2}^{[j]}(t)$; (b) $E_{H^1}^{[j]}(t)$, $j=0,1,2$.}
	\label{fig:8}
\end{figure}

\subsubsection{Example 3: The nonlocal-in-time fractional wave problem in 2D heterogeneous structures with  layered inclusions}
\label{subsec:6.3}

The third medium contains horizontal inclusion layers in the square domain $\Omega=[0,1\,\mathrm{cm}]^2$, with $\eps=1/16$. Each PUC contains a centered inclusion strip extending between its opposite vertical faces, with matrix material above and below. The physical inclusion layer has thickness $0.01875\,\mathrm{cm}$. Figure~\ref{fig:9} shows the composite and its PUC. Table~\ref{tab:materials} lists the material parameters for Example~3. The wave response is simulated over $t\in[0,1]\,\mathrm{s}$ using the time step $\tau=0.005\,\mathrm{s}$.

\begin{figure}[!htbp]
	\centering
	\includegraphics[width=0.82\textwidth]{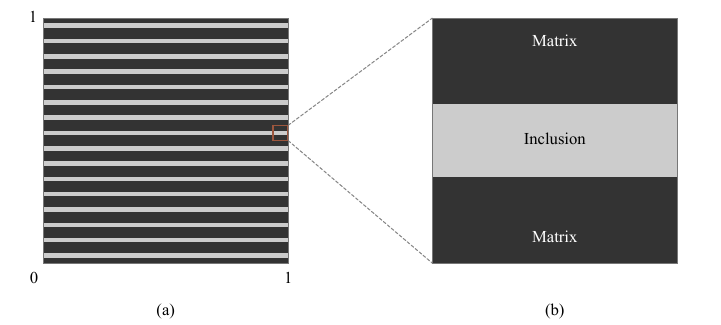}
	\caption{Example~3. (a) The 2D composite material $\Omega$; (b) reference PUC $Y$.}
	\label{fig:9}
\end{figure}

The first- and second-order auxiliary cell problems are solved with periodic boundary conditions. Boundary-layer corrections within the outermost layer of cells ensure that the reconstructed LOMS and HOMS solutions satisfy homogeneous Dirichlet conditions on the physical outer boundary.

The source term is
\begin{equation*}
f(\bx,t)=200\,\mathrm{kg}/(\mathrm{cm}\cdot\mathrm{s}^{2}),
\qquad \bx\in\Omega,\quad 0<t\le T.
\end{equation*}
Table~\ref{tab:cost} summarizes the mesh statistics and CPU times for Example~3. Figure~\ref{fig:10} compares the homogenized, LOMS, HOMS, and DNS solutions at $t_{36+\sigma}=0.18375\,\mathrm{s}$, while Figure~\ref{fig:11} shows the corresponding absolute errors. Figure~\ref{fig:12} presents the time evolution of the relative $L^2$-errors and relative $H^1$-seminorm errors.

The HOMS solution more closely reproduces the layer-dependent variations of the DNS displacement field, and the second-order correction reduces the displacement errors. Over most of the displayed time interval, the HOMS solution gives lower relative errors than the lower-order approximate solutions, particularly in the $H^1$ seminorm. The DNS mesh contains $256$ times as many elements as the homogenized mesh, substantially reducing the size of the macroscopic system advanced in time. The HOMS computation achieves a single-run speedup of approximately $2.53$ over DNS. For the reported mesh and time step, these results support the accuracy and efficiency of the HOMS asymptotic solution for this layered composite under constant loading.

\begin{figure}[!htbp]
	\centering
	\includegraphics[width=0.92\textwidth]{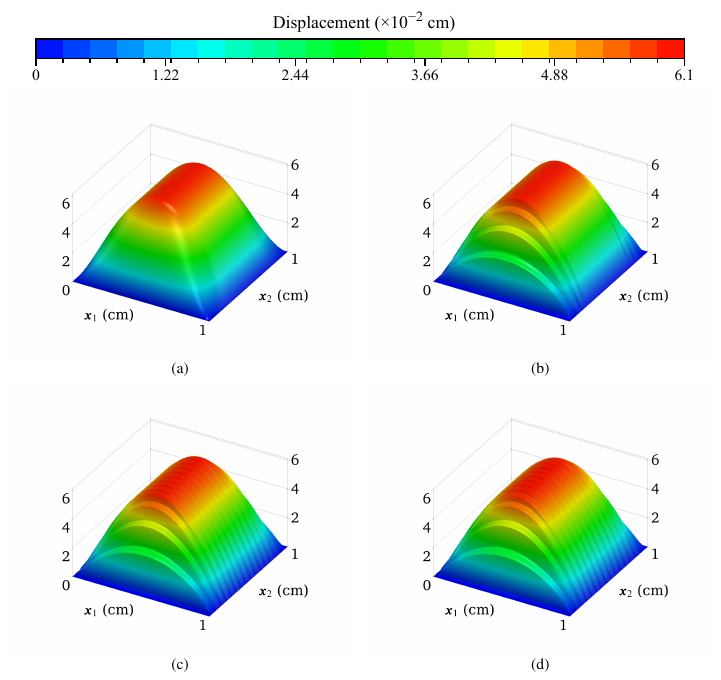}
	\caption{Example~3. Displacement fields: (a) $u_0$; (b) $u^{(1\varepsilon)}$; (c) $u^{(2\varepsilon)}$; (d) $u_{\rm DNS}^{\varepsilon}$.}
	\label{fig:10}

    \vspace{6pt}
\includegraphics[width=0.92\textwidth]{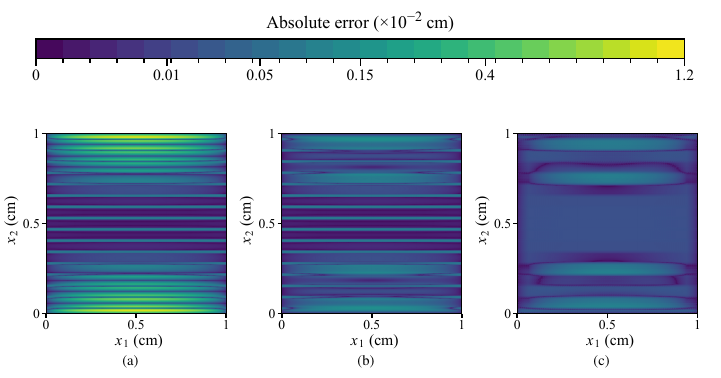}
\caption{Example~3. Absolute errors of displacement fields:
	(a) $|u_0-u_{\rm DNS}^{\varepsilon}|$;
	(b) $|u^{(1\varepsilon)}-u_{\rm DNS}^{\varepsilon}|$;
	(c) $|u^{(2\varepsilon)}-u_{\rm DNS}^{\varepsilon}|$.}
	\label{fig:11}
\end{figure}

\begin{figure}[!htbp]
	\centering
	\includegraphics[width=0.96\textwidth]{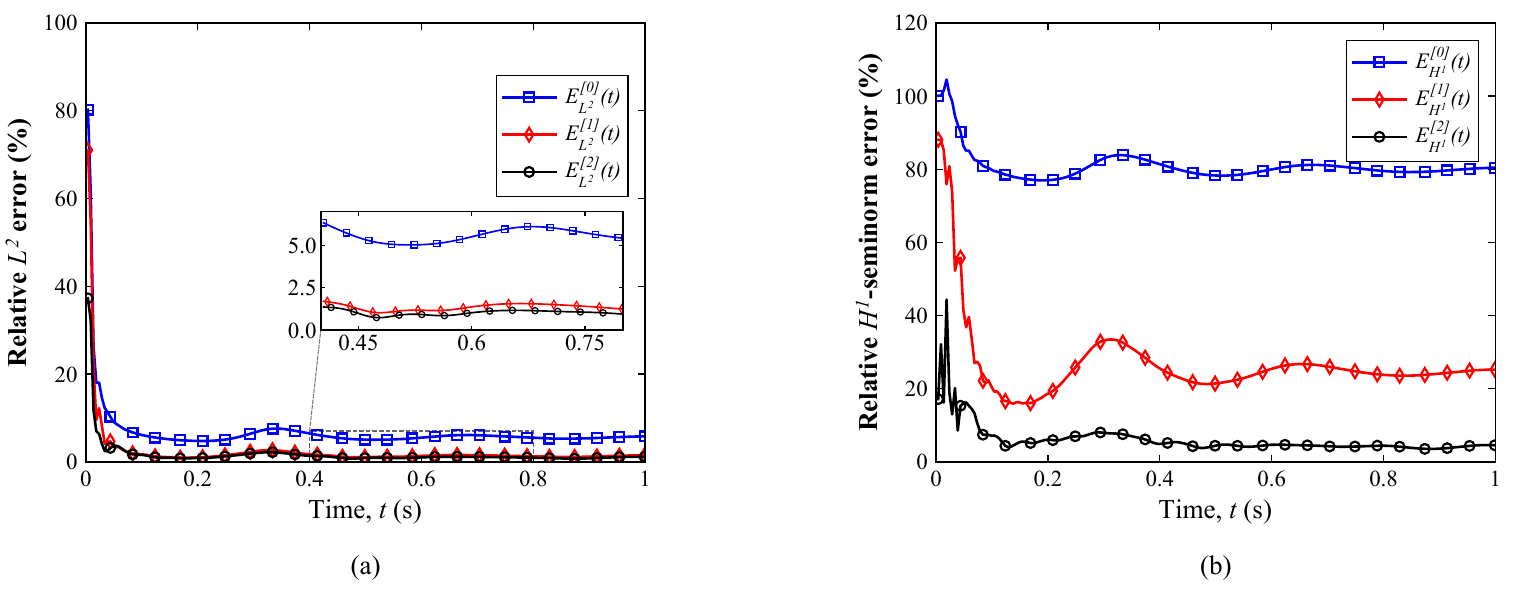}
	\caption{Example~3. Time evolution of relative errors: (a) $E_{L^2}^{[j]}(t)$; (b) $E_{H^1}^{[j]}(t)$, $j=0,1,2$.}
	\label{fig:12}
\end{figure}

\subsubsection{Example 4: The nonlocal-in-time fractional wave problem in 3D heterogeneous structures with cubic inclusions}
\label{subsec:6.4}

The fourth medium contains cubic inclusions in the domain $\Omega=[0,1\,\mathrm{cm}]^3$, with $\eps=1/6$. In the reference PUC $Y=[0,1]^3$, the centered cubic inclusion occupies $[1/4,3/4]^3$, and the surrounding region is the matrix. Figure~\ref{fig:13} shows the composite and its PUC. Table~\ref{tab:materials} lists the material parameters for Example~4. The wave response is simulated over $t\in[0,1]\,\mathrm{s}$ using the time step $\tau=0.005\,\mathrm{s}$.

\begin{figure}[!htbp]
	\centering
	\includegraphics[width=0.82\textwidth]{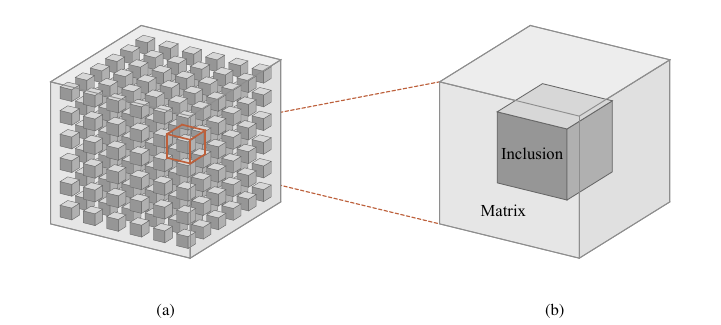}
	\caption{Example~4. (a) The 3D composite material $\Omega$; (b) reference PUC $Y$.}
	\label{fig:13}
\end{figure}

The source term is
\[
f(\bx,t)=600\,\mathrm{kg}/(\mathrm{cm}^{2}\cdot\mathrm{s}^{2}),
\qquad \bx\in\Omega,\quad 0<t\le T.
\]
Table~\ref{tab:cost} summarizes the mesh statistics and CPU times for Example~4. Figure~\ref{fig:14} compares the homogenized, LOMS, HOMS, and DNS solutions at $t_{50+\sigma}=0.25375\,\mathrm{s}$ on the horizontal section $x_3=7/12\,\mathrm{cm}$. Figure~\ref{fig:15} shows the corresponding absolute errors. Figure~\ref{fig:16} presents the time evolution of the relative $L^2$-errors and relative $H^1$-seminorm errors. Figure~\ref{fig:17} shows the displacement fields on four parallel sections at $x_1=1/4,5/12,7/12,$ and $3/4\,\mathrm{cm}$.

Figure~\ref{fig:16} shows that the HOMS solution has the smallest relative $L^2$-errors and relative $H^1$-seminorm errors after the initial peak.  The DNS mesh contains approximately $110.59$ times as many elements as the homogenized mesh. The smaller macroscopic system reduces the spatial cost of advancing the wave field and its fractional history at each time step. The HOMS computation achieves a single-run speedup of approximately $4.94$ over DNS. These results demonstrate the computational benefit of the second-order multi-scale asymptotic solution at the reported resolution.

\begin{figure}[!htbp]
	\centering
	\includegraphics[width=0.92\textwidth]{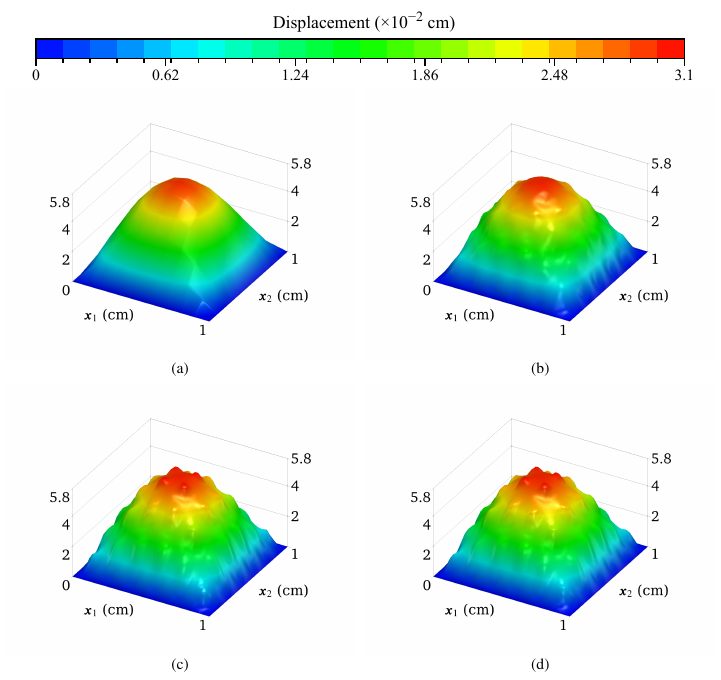}
	\caption{Example~4. Displacement fields: (a) $u_0$; (b) $u^{(1\varepsilon)}$; (c) $u^{(2\varepsilon)}$; (d) $u_{\rm DNS}^{\varepsilon}$.}
	\label{fig:14}

    \vspace{6pt}
\includegraphics[width=0.92\textwidth]{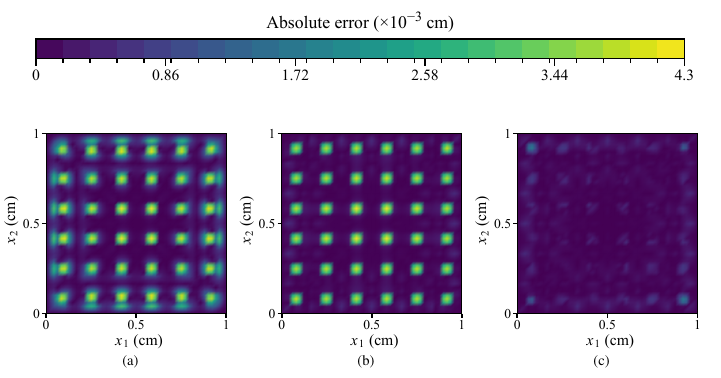}
\caption{Example~4. Absolute errors of displacement fields:
	(a) $|u_0-u_{\rm DNS}^{\varepsilon}|$;
	(b) $|u^{(1\varepsilon)}-u_{\rm DNS}^{\varepsilon}|$;
	(c) $|u^{(2\varepsilon)}-u_{\rm DNS}^{\varepsilon}|$.}
	\label{fig:15}
\end{figure}

\begin{figure}[!htbp]
	\centering
	\includegraphics[width=0.96\textwidth]{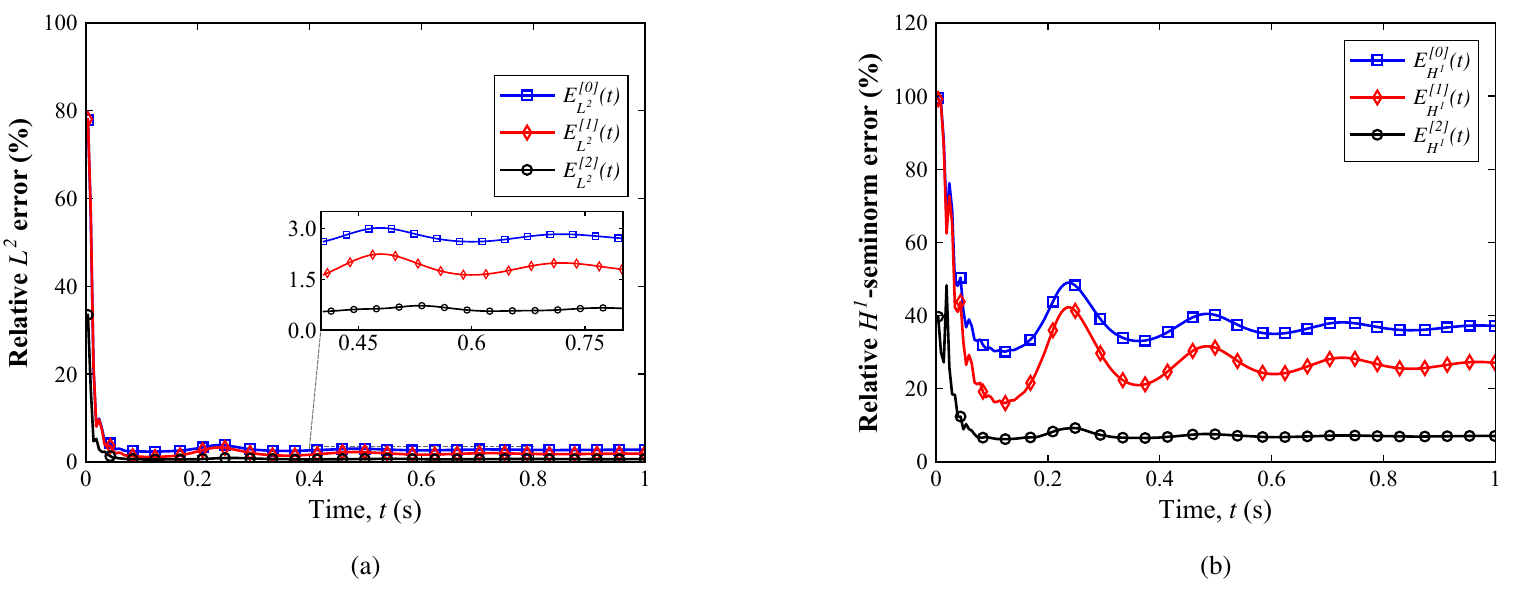}
	\caption{Example~4. Time evolution of relative errors: (a) $E_{L^2}^{[j]}(t)$; (b) $E_{H^1}^{[j]}(t)$, $j=0,1,2$.}
	\label{fig:16}

    \vspace{6pt}
\includegraphics[width=0.94\textwidth]{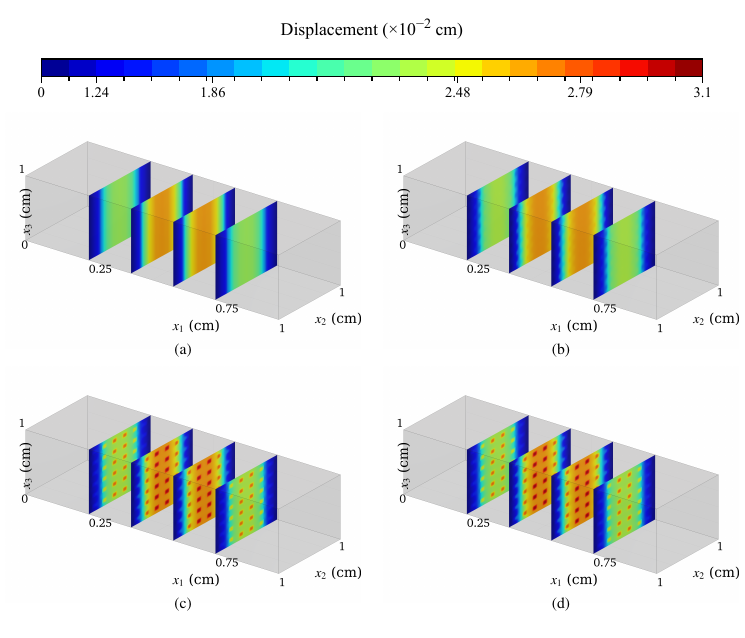}
	\caption{Example~4. Displacement fields on four vertical sections: (a) $u_0$; (b) $u^{(1\varepsilon)}$; (c) $u^{(2\varepsilon)}$; (d) $u_{\rm DNS}^{\varepsilon}$.}
	\label{fig:17}
\end{figure}

Across Examples~1--4, the homogenized solution captures the smooth macroscopic response, while the LOMS solution recovers part of the microscale variation. The HOMS solution more closely reproduces the DNS displacement fields in Figures~\ref{fig:2}, \ref{fig:6}, \ref{fig:10}, and~\ref{fig:14}.

\subsection{Applications of the HOMS computational method}
\subsubsection{Example 5: Simulation of a nonlocal-in-time fractional wave problem in a 3D heterogeneous composite plate with a large number of cubic inclusions}
\label{subsec:6.5}

The fifth medium contains cubic inclusions in the plate domain $\Omega=[0,10\,\mathrm{cm}]\times[0,10\,\mathrm{cm}]\times[0,2\,\mathrm{cm}]$, with $\eps=1/4$. The cubic-inclusion geometry and material parameters are the same as in Example~4. In each reference PUC, the centered inclusion occupies $[1/4,3/4]^3$. Figure~\ref{fig:18} shows the composite and its PUC. The simulation covers $0\le t\le1\,\mathrm{s}$ with $\tau=0.005\,\mathrm{s}$.

\begin{figure}[!htbp]
	\centering
	\includegraphics[width=0.82\textwidth]{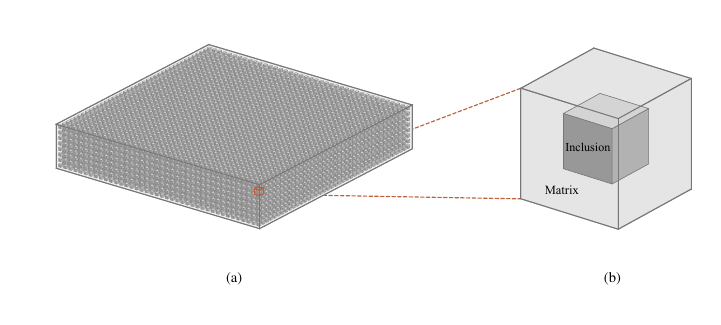}
	\caption{Example~5. (a) The 3D heterogeneous composite plate $\Omega$; (b) reference PUC $Y$.}
	\label{fig:18}
\end{figure}

The source term is the same as in Example~4. DNS would require a fine mesh throughout the large domain to resolve the repeated material interfaces. The estimated memory and computational requirements of the DNS mesh exceed the available resources, preventing a converged reference solution from being obtained. Table~\ref{tab:cost} reports the estimated DNS mesh sizes and the cell and homogenized mesh sizes for Example~5.

Figure~\ref{fig:19} compares the three solutions on the horizontal section $x_3=9/8\,\mathrm{cm}$ at $t_{50+\sigma}=0.25375\,\mathrm{s}$, together with enlarged views of the interior region $4\le x_1,x_2\le6\,\mathrm{cm}$. To examine the response through the thickness, Figure~\ref{fig:20} shows the three solutions on the vertical section $x_2=41/8\,\mathrm{cm}$ at the same time. The profiles along $L_1$--$L_4$, located at $x_1=17/8,33/8,49/8,65/8\,\mathrm{cm}$, compare their local variations along $x_3$.

\begin{figure}[!htbp]
	\centering
	\includegraphics[width=\textwidth]{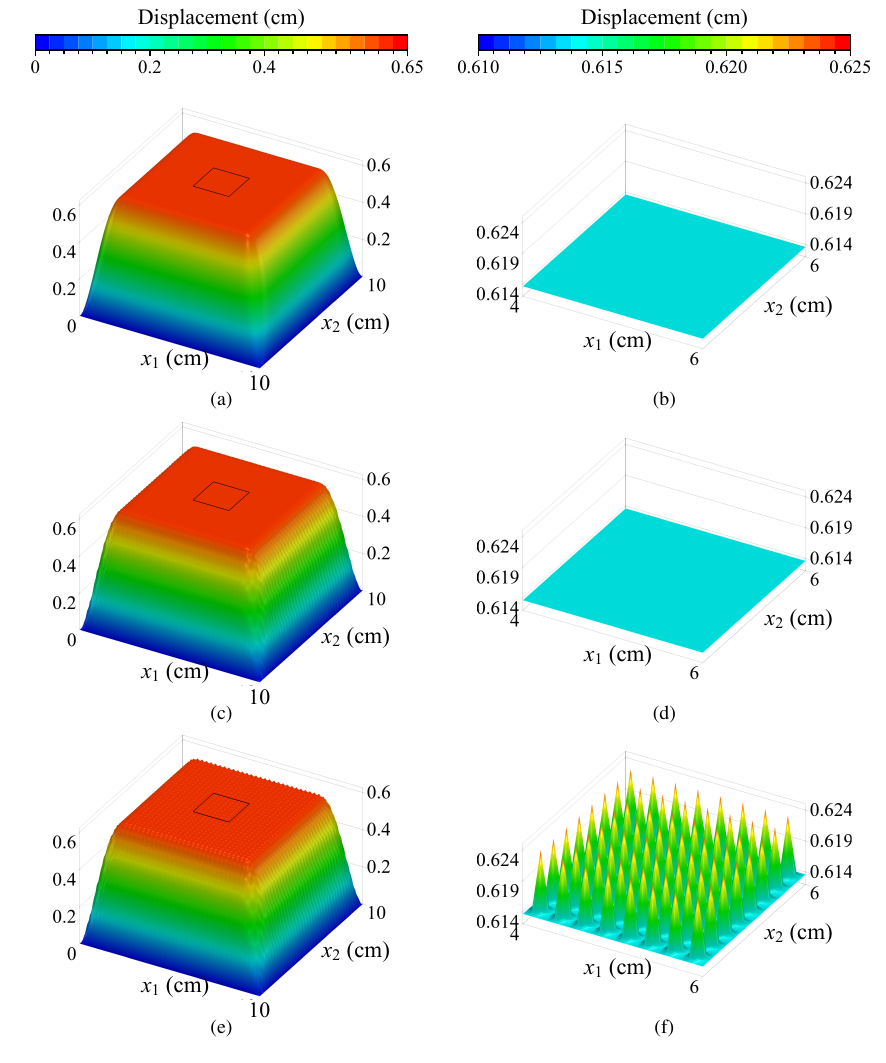}
	\caption{Example~5. Displacement fields on the horizontal section: (a,b) $u_0$; (c,d) $u^{(1\varepsilon)}$; (e,f) $u^{(2\varepsilon)}$. Panels (b,d,f) show enlarged views.}
	\label{fig:19}
\end{figure}

\begin{figure}[!htbp]
	\centering
	\includegraphics[width=\textwidth]{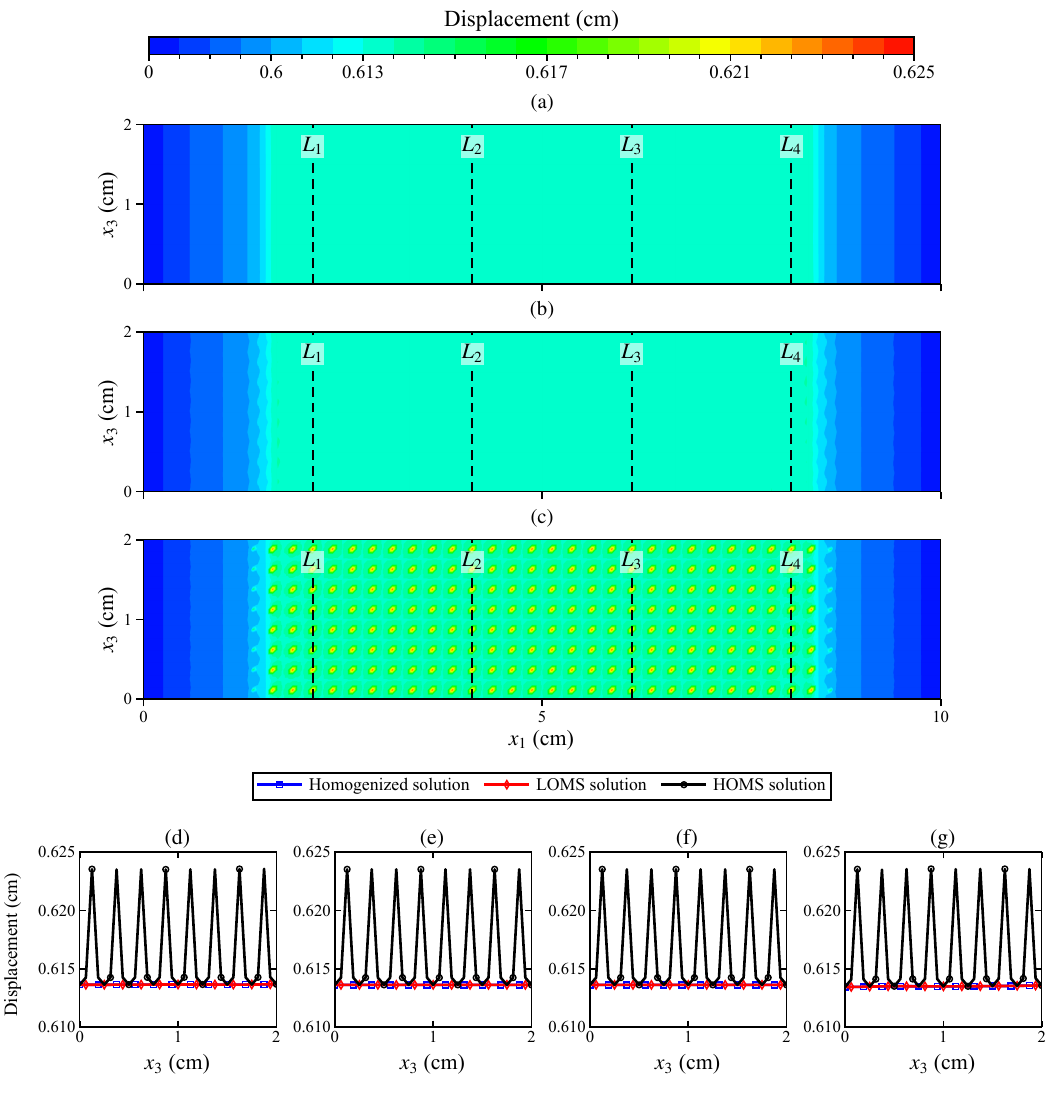}
	\caption{Example~5. Displacement fields and line profiles: (a) $u_0$; (b) $u^{(1\varepsilon)}$; (c) $u^{(2\varepsilon)}$; (d)--(g) profiles along $L_1$--$L_4$, respectively.}
	\label{fig:20}
\end{figure}

\endgroup
\section{Conclusions}
\label{sec:7}

This paper presents a computationally efficient HOMS method for nonlocal-in-time wave propagation in composite materials with Caputo fractional damping. The method combines the macroscopic homogenized equation with first- and second-order auxiliary cell problems to construct multi-scale asymptotic solutions for microscopic wave fields. Moreover, local pointwise error analysis shows that the HOMS solution preserves local physical balance and cancels the leading residual of the LOMS solution. The global error estimates also establish $O(\varepsilon)$-order convergence in the energy norm. Furthermore, a novel fully discrete $L2$-$1_\sigma$--Newmark scheme is developed for the homogenized equation, with unconditional shifted-energy stability under the stated assumptions. The fully discrete homogenized displacement and velocity converge to their continuous counterparts with errors of $O(h_Y^2+h^2+\tau^2)$.

Two- and three-dimensional numerical examples demonstrate the accuracy of the HOMS solution and the computational efficiency of the HOMS method. The HOMS solution gives smaller relative errors than the homogenized and LOMS solutions over most of the displayed interval. Based on the total computational cost of the HOMS method, reported speedups over DNS range from $2.53$ to $9.51$. Moreover, the numerical example for composite structures with massive unit cells demonstrates the computational feasibility and robustness of the HOMS method for large-scale engineering simulation.

The current error estimates rely on high temporal regularity assumptions. Future work will investigate corresponding error estimates under weaker temporal regularity assumptions. We will also explore extending the HOMS method to simulations of materials described by fractional viscoelastic models and coupled multiphysics problems involving fractional derivatives.

\section*{ Declarations }

\textbf{Author Contributions} Gehao Wang: Conceptualization, Methodology, Software, Validation, Formal analysis, Writing–original draft, Writing–review \& editing, Visualization.
Hao Dong: Conceptualization, Methodology, Formal analysis, Writing–original draft, Writing–review \& editing, Supervision, Project administration, Funding acquisition.
Yufeng Nie: Methodology, Formal analysis, Writing–review \& editing.
Junzhi Cui: Conceptualization, Methodology, Supervision, Writing–review \& editing.

\textbf{Funding}  This research was supported by the National Natural Science Foundation of China (No. 12471387), Natural Science Basic Research Program of Shaanxi (No. 2026JC-YXQN-017), Xidian University Specially Funded Project for Interdisciplinary Exploration (No. TZJH2024008), Fundamental Research Funds for the Central Universities (No. QTZX25082), Innovation Capability Support Program of Shaanxi Province (No. 2024RS-CXTD-88).

\textbf{Data Availability} The data that support this study are available from the corresponding author on reasonable request.

\textbf{Conflict of interest} The authors declare they have no conflicts of interest.


%
%



\end{document}